\documentstyle[ctagsplt,amscd]{amsart}

\numberwithin{equation}{section}

\newtheorem{theorem}{Theorem}[section]
\newtheorem{lemma}[theorem]{Lemma}
\newtheorem{proposition}[theorem]{Proposition}
\newtheorem{corollary}[theorem]{Corollary}

\newtheorem{claim}[theorem]{Claim}

\theoremstyle{definition}
\newtheorem{definition}[theorem]{Definition}

\theoremstyle{remark}
\newtheorem{remark}[theorem]{Remark}
\newtheorem{example}[theorem]{Example}
\newtheorem{acknowledgement}{Acknowledgement}

\renewcommand\Re{\operatorname{Re}}
\renewcommand\Im{\operatorname{Im}}
\newcommand\grad{\operatorname{grad}}

\newcommand\vol{\operatorname{vol}}
\newcommand\ch{\operatorname{ch}}

\newcommand\Ad{\operatorname{Ad}}
\newcommand\GL{\operatorname{GL}} %
\newcommand\SL{\operatorname{SL}} 
\newcommand\U{\operatorname{U}} %
\newcommand\id{\operatorname{id}}
\newcommand\Hom{\operatorname{Hom}} 
\newcommand\Aut{\operatorname{Aut}} 

\newcommand\eps{\varepsilon}

\newcommand\bu{{\mbox{\large$.$}}}

\newcommand\co{^{\bold C}} 
\newcommand\st{^{\mathrm s}} 
\newcommand\sst{^{\mathrm s\mathrm s}} 
\newcommand\qu{/\kern-.7ex/} 
\newcommand\sq{\sqrt{-1}}

\newcommand\romquote{\rom`\kern-.2ex\rom`\kern-.2ex}
\newcommand\romunquote{\rom'\kern-.2ex\rom'} 

\begin{document}

\title[Slices, Reduction and Multiplicities]{Holomorphic Slices,
Symplectic Reduction and Multiplicities of Representations}
\author{Reyer Sjamaar}
\address{University of Pennsylvania, Department of Mathematics,
Philadelphia, Pennsylvania 19104}
\curraddr{Massachusetts Institute of Technology, Department of
Mathematics, Cambridge, Massachusetts 02139-4307}
\email{sjamaar@@math.mit.edu}
\thanks{This work was partially supported by NSF grant DMS92-03398.}
\date{March 1993}

\maketitle

\begin{abstract}
I prove the existence of slices for an action of a reductive complex
Lie group on a K\"ahler manifold at certain orbits, namely those
orbits that intersect the zero level set of a momentum map for the
action of a compact real form of the group. I give applications of
this result to symplectic reduction and geometric quantization at
singular levels of the momentum map. In particular, I obtain a
formula for the multiplicities of the irreducible representations
occurring in the quantization in terms of symplectic invariants of
reduced spaces, generalizing a result of Guillemin and Sternberg.
\end{abstract}

\section*{Introduction}

There has recently been much interest in formulas for multiplicities
of Lie group representations arising in various different ways from
group actions on manifolds. Typically, one can think of the manifold
as being the phase space $M$ of a classical physical system acted upon
by a group $G$ of symmetries, from which one obtains a unitary
representation of $G$ through some kind of ``quantization''. A
prototype of such formulas is the multiplicity formula of Guillemin
and Sternberg \cite{gu:ge}. In their set-up the space $M$ is a compact
K\"ahler manifold on which the compact group $G$ acts by holomorphic
transformations, and the associated representation of $G$ is the space
of holomorphic sections of a certain $G$-equivariant line bundle over
$M$ (``geometric quantization''). The main result of \cite{gu:ge}
expresses the multiplicities of the irreducible components of this
representation in terms of the Riemann-Roch numbers of the symplectic
(or Marsden-Weinstein-Meyer) quotients $M_\lambda$ of $M$. The most
important auxiliary result is that the symplectic quotient at the zero
level, $M_0$, can be identified with a geometric quotient of $M$ by
the complexified group $G\co$ as defined by Mumford \cite{mu:ge}.

The purpose of this paper is to generalize these results to the case
where the symplectic quotient $M_\lambda$ is singular, a case which is
of some interest in applications, but was excluded by Guillemin and
Sternberg. This involves a closer study of the orbit structure of the
action of the reductive group $G\co$ on $M$. The main technical
result, expounded in Section \ref{section:slices}, is that one can
construct slices for the $G\co$-action at points that are in the zero
level set of a momentum map. The proof of this holomorphic slice
theorem utilizes H\"ormander's $L_2$-estimates for the Cauchy-Riemann
operator. For affine algebraic manifolds it is a special case of
Luna's \'etale slice theorem \cite{lu:sl}.

Kirwan \cite{ki:coh} has introduced the notion of the K\"ahler
quotient of $M$ by $G\co$, which is the K\"ahler analogue of Mumford's
categorical quotient. It is roughly speaking defined as the space of
closed $G\co$-orbits in $M$. Kirwan showed it is homeomorphic to the
symplectic quotient $M_0$, generalizing the result of Guillemin and
Sternberg referred to above. In \cite{sj:st} it was shown that $M_0$
is a so-called symplectic stratified space. In Section
\ref{section:quotient} I exploit the holomorphic slice theorem to
study the analytic structure of Kirwan's quotient and to compare it to
the stratified symplectic structure of the symplectic quotient.

In the case of a projective manifold $M$ endowed with the Fubini-Study
symplectic form, Kirwan and, independently, Ness \cite{ne:st} showed
that the symplectic quotient $M_0$ coincides with Mumford's
categorical quotient, the Proj of the invariant part of the
homogeneous coordinate ring of $M$. I show that the same conclusion
holds when $M$ has an arbitrary integral K\"ahler structure. (Under
this hypothesis $M$ has a unique algebraic structure by Kodaira's
Embedding Theorem, but the symplectic structure is not necessarily one
coming from a Fubini-Study metric.) This result allows me to carry
through the geometric quantization of the symplectic quotient. Now a
theorem of Boutot \cite{bo:si} asserts that the singularities of a
quotient such as $M_0$ are {\em rational}. This basically says that
the Riemann-Roch numbers of $M_0$ are equal to those of any blowup,
which, finally, leads to a generalization of Guillemin and Sternberg's
multiplicity formula.

\begin{acknowledgement}
I am grateful to Hans Duistermaat for suggesting to me the problem
discussed in this paper and for showing me an unpublished manuscript
of his, a succinct version of which appeared in \cite{du:mu}. Eugene
Lerman has been a great help in carrying out this work. Part of it
appears as a joint announcement in \cite{le:re}. I would also like to
thank Eugenio Calabi and Charlie Epstein for their generous help with
the material in Section \ref{subsection:totallyreal}. I have
furthermore benefited from helpful discussions with Victor Guillemin
and Viktor Ginzburg.
\end{acknowledgement}

\section{Holomorphic Slices}\label{section:slices}

Let $X$ be a complex space and let $G\co$ be a reductive complex Lie
group acting holomorpically on $X$. I think of $G\co$ as being the
complexification of a compact real Lie group $G$.
\begin{definition}
A {\em slice\/} at $x$ for the $G\co$-action is a locally closed
analytic subspace $S$ of $X$ with the following properties:
\begin{enumerate}
\item $x\in S$;
\item $G\co S$ of $S$ is open in $X$;
\item $S$ is invariant under the action of the stabilizer $(G\co)_x$;
\item\label{bundle} the natural $G\co$-equivariant map from the
associated bundle $G\co\times_{(G\co)_x}S$ into $X$, which sends an
equivalence class $[g,y]$ to the point $gy$, is an analytic
isomorphism onto $G\co S$.
\end{enumerate}
\end{definition}
It follows from (\ref{bundle}) that for all $y$ in $S$ the stabilizer
$(G\co)_y$ is contained in $(G\co)_x$. Furthermore, if $X$ is
nonsingular at $x$, a slice $S$, if it exists, is nonsingular at $x$
and transverse to the orbit $G\co x$.

The problem of constructing slices has been solved by Luna
\cite{lu:sl} for affine varieties and by Snow \cite{sn:re} for Stein
spaces. One difficulty of the problem lies in the fact that an action
of $G\co$ is typically not proper, unless it is locally free. One
therefore faces the challenge of controlling the behaviour of the
action ``at infinity in the group''. Another snag is that there may be
cohomological obstructions to analytically embedding the ``normal
bundle'' $G\co\times_{(G\co)_x}S$ of the orbit $G\co x$ into $X$.
These obstructions vanish if the orbit is (analytically isomorphic to)
an affine variety. A theorem of Matsushima's \cite{ma:es} says that a
$G\co$-orbit is affine if and only if the isotropy subgroup $(G\co)_x$
is reductive. But even if the isotropy of $x$ is reductive one cannot
always construct a slice at $x$. (Cf.\ Richardson's example
\cite{le:re,lu:sl,ri:de} of the standard action of $\SL(2)$ on
homogeneous cubic polynomials, and also Trautman \cite{tr:or}.) The
additional condition that Luna, resp.\ Snow, impose in the context of
an affine variety $X$, resp. a Stein space $X$, in order to deduce the
existence of a slice is that the orbit should be {\em closed\/} in
$X$.

The above notion of a slice is slightly weaker than that of Luna and
Snow, who require the set $G\co S$ to be saturated with respect to a
quotient mapping. In our context the definition of a quotient depends
upon the choice of a momentum map. In the next section we shall see
that for any choice of a momentum map there always exists slices $S$
such that $G\co S$ is saturated with respect to the corresponding
quotient map (Proposition \ref{proposition:saturated}).

In this section I demonstrate the existence of slices at certain
affine orbits of a $G\co$-action on a K\"ahler manifold. I was led to
this result by the striking resemblance between Luna's and Snow's
slice theorems and the normal forms in symplectic geometry due to
Marle \cite{ma:mo} and Guillemin and Sternberg \cite{gu:no}. Before
formulating the theorem I have to state a number of definitions and
auxiliary results. In Section \ref{subsection:orbits} I discuss
momentum maps on K\"ahler manifolds and the notion of orbital
convexity. Section \ref{subsection:totallyreal} contains a result
concerning interpolation between K\"ahler metrics in the neighbourhood
of a totally real submanifold of a complex manifold, which relies on
H\"ormander's $\bar\partial$-estimates, and which is the main
ingredient in the proof of the slice theorem. In Section
\ref{subsection:slices} I prove the slice theorem and discuss some of
its immediate consequences.

\subsection{Orbital convexity and isotropic
orbits}\label{subsection:orbits}

Recall that the decomposition of the complexified Lie algebra $\frak
g\co=\frak g\otimes\bold C$ into a direct sum $\frak g\co=\frak
g\oplus\sq\,\frak g$ gives rise to the polar (or Cartan) decomposition
$G\co=G\cdot\exp\sq\,\frak g$. The map $G\times\sq\,\frak g\to G\co$
sending $(k,\sq\,\xi)$ to $k\exp\sq\,\xi$ is a diffeomorphism onto,
and every element $g$ of $G\co$ can be uniquely decomposed into a
product $g=k\exp\sq\,\xi$, with $k\in G$ and $\xi\in\frak g$.

\begin{definition}[Heinzner {\cite{he:ge}}]
A subset $A$ of a $G\co$-space $X$ is called {\em orbitally convex\/}
with respect to the $G\co$-action if it is $G$-invariant and if for
all $x$ in $U$ and all $\xi$ in $\frak g$ the intersection of the
curve $\bigl\{\,\exp(\sq\,t\xi)x:t\in\bold R\,\bigr\}$ with $A$ is
connected.  Equivalently, $A$ is orbitally convex if and only if it is
$G$-invariant and for all $x$ in $A$ and all $\xi$ in $\frak g$ the
fact that both $x$ and $\exp(\sq\,\xi)x$ are in $A$ implies that
$\exp(\sq\,t\xi)x\in A$ for all $t\in[0,1]$.
\end{definition}

\begin{remark}\label{remark:trivial}
If $f\colon X\to Y$ is a $G\co$-equivariant map between $G\co$-spaces
$X$ and $Y$, and $C$ is an orbitally convex subset of $Y$, then it
follows immediately from the definition that $f^{-1}(C)$ is orbitally
convex in $X$.
\end{remark}

A $G$-equivariant map defined on an orbitally convex open set can be
analytically continued to a $G\co$-equivariant map.

\begin{proposition}[Heinzner {\cite{he:ge}}, Koras
\cite{ko:li}]\label{proposition:orbit}
Let $X$ and $Y$ be complex manifolds acted upon by $G\co$. If $A$ is
an orbitally convex open subset of $X$ and $f\colon A\to Y$ is a
$G$-equivariant holomorphic map\rom, then $f$ can be uniquely extended
to a $G\co$-equivariant holomorphic map $f\co\colon G\co A\to Y$.

Consequently, if the image $f(A)$ is open and orbitally convex in $Y$
and $f\colon A\to f(A)$ is biholomorphic\rom, then the extension
$f\co\colon G\co A\to Y$ is biholomorphic onto the open subset $G\co
f(A)$.
\end{proposition}

\begin{pf}
The only way to extend $f$ equivariantly is by putting
$f\co\bigl(g\exp(\sq\,\xi)x\bigr)=g\exp(\sq\,\xi)f(x)$ for all $x$ in
$A$, $g$ in $G$ and $\xi$ in $\frak g$. We have to check this is
well-defined.

Let $x\in A$ and $\xi\in\frak g$ be such that $\exp(\sq\,\xi)x\in A$.
Then by assumption $\exp(\sq\,t\xi)x\in A$ for all $t$ between 0 and
1. So $f\bigl(\exp(\sq\,t\xi)x\bigr)$ is well-defined for $0\leq t\leq
1$.  Define the curves $\alpha(t)$ and $\beta(t)$ in $Y$ by
$\alpha(t)=f\bigl(\exp(\sq\,t\xi)x\bigr)$ and $\beta(t)=
\exp(\sq\,t\xi)f(x)$ for $0\leq t\leq 1$.  Then $\alpha(t)$ and
$\beta(t)$ are integral curves of the vector fields $f_*(\sq\,\xi)_X$
and $(\sq\,\xi)_Y$ respectively, both with the same initial value
$f(x)$. Now since $f$ is $G$-equivariant we have $f_*\xi_X=\xi_Y$,
and, because $f$ is also holomorphic, $f_*(\sq\,\xi)_X= f_*(J\xi_X)=
Jf_*\xi_X= J\xi_Y= (\sq\,\xi)_Y$. Hence $\alpha(t)=\beta(t)$, in other
words $f\bigl(\exp(\sq\,t\xi)x\bigr)=\exp(\sq\,t\xi)f(x)$ for $0\leq
t\leq 1$.

It follows that for all $x$ in $A$ and all $\xi$ in $\frak g$ such
that $\exp(\sq\,\xi)x$ is in $A$ we have
$f\bigl(\exp(\sq\,\xi)x\bigr)=\exp(\sq\,\xi)f(x)$. It is easy to
deduce from this that $f\co$ is well-defined.

Finally observe that if the image $f(A)$ is open and orbitally convex
in $Y$ and $f\colon A\to f(A)$ is biholomorphic, then the inverse
$f^{-1}$ also has a holomorphic extension $(f^{-1})\co\colon G\co
f(A)\to G\co A$, and by uniqueness this must be the inverse of $f\co$.
\end{pf}

\begin{remark}\label{remark:local}
Suppose we drop the assumption that $A$ is orbitally convex from the
statement of the proposition. Then it is not necessarily true that
$f\bigl(\exp(\sq\,t\xi)x\bigr)$ is equal to $\exp(\sq\,t\xi)f(x)$ for
all $t$ such that $\exp(\sq\,t\xi)x\in A$. But if we put
$I=\{\,t\in\bold R:\exp(\sq\,t\xi)x\in A\,\}$ and let $I^0$ be the
connected component of $I$ containing $0$, then the above proof shows
that $f\bigl(\exp(\sq\,t\xi)x\bigr)=\exp(\sq\,t\xi)f(x)$ for all $t$
in $I^0$.
\end{remark}

In the remainder of this section $M$ shall denote a K\"ahler manifold,
not necessarily compact, with infinitely differentiable K\"ahler
metric $ds^2$, K\"ahler form $\omega=-\Im ds^2$, and complex structure
$J$. Then $\Re ds^2=\omega(\cdot,J\cdot)$ is the corresponding
Riemannian metric. We may assume without loss of generality that
$ds^2$ is invariant under the compact group $G$. So the
transformations on $M$ defined by $G$ are holomorphic and they are
isometries with respect to the K\"ahler metric. The action of $G$ is
called {\em Hamiltonian\/} if for all $\xi$ in the Lie algebra $\frak
g$ of $G$ the vector field $\xi_M$ on $M$ induced by $\xi$ is
Hamiltonian. In this case we have a {\em momentum map\/} $\Phi$ from
$M$ to the dual $\frak g^*$ of the Lie algebra of $G$ with the
property that
$$d\Phi^\xi=\iota_{\xi_M}\omega$$
for all $\xi$. Here $\Phi^\xi$ is the $\xi$-th component of $\Phi$,
defined by $\Phi^\xi(m)=\bigl(\Phi(m)\bigr)(\xi)$. After averaging
with respect to the given action on $M$ and the coadjoint action on
$\frak g^*$ we may assume that the map $\Phi$ is $G$-equivariant. An
equivariant momentum map is uniquely determined up to additive
constants ranging over the $\Ad^*$-fixed vectors in $\frak g^*$. (So
if $G$ is connected, the number of degrees of freedom is equal to the
dimension of the centre of $G$.) It is easy to give sufficient
conditions for the existence of a momentum map, for instance, the
first Betti number of $M$ is zero, or the K\"ahler form $\omega$ is
exact.  (See e.g.\ \cite{gu:sy,we:le}.) More surprisingly, by a
theorem of Frankel \cite{fr:fi} a momentum map always exists if the
action has at least one fixed point and $M$ is {\em compact}. A
necessary and sufficient condition for a holomorphic $G$-action on a
compact K\"ahler manifold to be Hamiltonian is that for every vector
$\xi\in\frak g$ the holomorphic vector field $\xi_M$ should be killed
by every global holomorphic one-form $\alpha$ on $M$,
$\alpha(\xi_M)=0$.  (This follows from Frankel's theorem and a fixed
point theorem of Sommese \cite{so:ex}.)  Note that this condition is
independent of the K\"ahler structure.

If $M$ is $\bold C^n$ with the standard Hermitian structure
$dS^2=\sum_idz_i\otimes d\bar z_i$ and the standard symplectic form
$\Omega=\sq\big/2\,\sum_idz_i\wedge d\bar z_i$, then a momentum map
$\Phi_{\bold C^n}$ is given by the formula
\begin{equation}\label{equation:quadratic}
\Phi_{\bold C^n}^\xi(v)=1/2\,\Omega\bigl(\xi_{\bold C^n}(v),v\bigr),
\end{equation}
where $\xi _{\bold C^n}$ denotes the image of $\xi\in\frak g$ in the
Lie algebra $\frak s\frak p({\bold C^n},\Omega)$, and $v\in\bold
C^n$.

Because $G$ acts holomorphically on $M$, there is a natural way to
define an action of the complexified Lie algebra $\frak g^*$: For any
$\xi$ in $\frak g$ the vector field $(\sq\,\xi)_M$ induced by
$\sq\,\xi$ is equal to $J\xi_M$. It follows easily from the definition
of a momentum map that $J\xi_M$ is equal to the gradient vector field
(with respect to the Riemannian metric $\Re ds^2$) of the $\xi$-th
component of the momentum map,
\begin{equation}\label{equation:grad}
(\sq\,\xi)_M=J\xi_M=\grad\Phi^\xi.
\end{equation}
We will assume that these vector fields are {\em complete\/} for all
$\xi$ in $\frak g$. This assumption implies that the action of $G$
extends uniquely to a holomorphic action of $G\co$. The assumption
holds for instance if $M$ is compact, or if $M$ is the total space of
a vector bundle over a compact manifold on which $G$ acts by vector
bundle transformations. The identity (\ref{equation:grad}) will enable
us to gain control over the behaviour of the action ``at infinity in
the group''. For one thing, it implies that the trajectory
$\gamma(t)$ of $\grad\Phi^\xi$ through a point $x$ in $M$ is given by
$\gamma(t)=\exp(\sq\,t\xi)x$, which does not depend on the choice of
the K\"ahler metric or the momentum map.

Here is another application of (\ref{equation:grad}). A submanifold
$X$ of $M$ is called {\em totally real\/} if $T_xX\cap
J\bigl(T_xX\bigr)=\{0\}$ for all $x\in X$.

\begin{proposition}\label{proposition:totallyreal}
Assume $G$ is connected\rom. Consider the following conditions on a
point $m\in M$\rom:
\begin{enumerate}
\item\label{fixed} $\Phi(m)$ is fixed under the coadjoint action of
$G$ on $\frak g^*$\rom;
\item\label{isotropic} The orbit $Gm$ is isotropic with respect to the
K\"ahler form\rom;
\item\label{stabilizer} The complex stabilizer $(G\co)_m$ of $m$ is
equal to the complexification $(G_m)\co$ of the compact stabilizer
$G_m$\rom, $(G_m)\co=(G\co)_m$\rom;
\item\label{realorbit} The $G$-orbit through $m$ is totally real\rom.
\end{enumerate}
Conditions \rom(\ref{fixed}\rom) and \rom(\ref{isotropic}\rom) are
equivalent\rom.  Any one of these conditions implies
\rom(\ref{stabilizer}\rom)\rom; and \rom(\ref{stabilizer}\rom) implies
\rom(\ref{realorbit}\rom)\rom.
\end{proposition}

\begin{pf}
Put $\mu=\Phi(m)$. Let $G_\mu$ be the stabilizer of $\mu$ with respect
to the coadjoint action; it is well-known that $G_\mu$ is a connected
subgroup of $G$. Let $G\mu$ be the coadjoint orbit through
$\mu$. We regard $G\mu$ as a symplectic manifold with the
Kirillov-Kostant-Souriau symplectic form. Denote the tangent space
$T_m(G m)$ to the compact orbit by $\frak m$; then $\frak m$ is
isomorphic to $\frak g/\frak g_m$ as an $H$-module. Similarly, let
$\frak n$ denote the tangent space $T_m(G_\mu m)$ to the orbit $G_\mu
m$; then $\frak n\cong\frak g_\mu/\frak g_m$. We have a fibration
$$
G_\mu m\hookrightarrow G m\overset{\Phi}{\to}G\mu,
$$
which on the tangent level leads to an exact sequence of vector spaces
$$
0\to\frak n\to\frak m\overset{d\Phi}{\to}T_\mu(G\mu)\to 0.
$$
The restriction of the symplectic form $\omega$ to $\frak m$ is an
alternating bilinear (``presymplectic'') form. It follows from the
fact that $\Phi$ is a Poisson map that $d\Phi$ preserves the
presymplectic forms. Since $T_\mu(G\mu)$ is symplectic, $\frak n$
is exactly the nullspace of $\omega|_{\frak m}$. Therefore, $\frak m$
is isotropic if and only if $T_\mu(G\mu)=0$, that is, $G_\mu=G$,
in other words, $\mu$ is $G$-fixed. This shows that (\ref{fixed}) is
equivalent to (\ref{isotropic}).

We now prove (\ref{fixed}) implies (\ref{stabilizer}). It is easy
to see that for any point $m$ in $M$ the complex stabilizer $(G\co)_m$
contains the complexification $(G_m)\co$ of the compact stabilizer
$G_m$. Now suppose $\mu$ is fixed under the coadjoint action. Let
$g\exp\sq\,\xi$ be an arbitrary element of $\dim(G\co)_m$, where $g\in
G$ and $\xi\in\frak g$. We want to show that $g\in G_m$ and
$\xi\in\frak g_m$. (Cf.\ Kirwan \cite{ki:coh} for this part of the
argument.) By $G$-equivariance of the momentum map we have
$$
\Phi\bigl(\exp(\sq\,\xi)m\bigr)=
g^{-1}\Phi\bigl(g\exp(\sq\,\xi)m\bigr)= g^{-1}\Phi(m)=g^{-1}\mu=\mu,
$$
and therefore $\Phi^\xi\bigl(\exp(\sq\,\xi)m\bigr)=\Phi^\xi(m)$. By
(\ref{equation:grad}) the curve $\exp(\sq\,t\xi)m$ is the
gradient trajectory of the vector field $(\sq\,\xi)_M$ through $m$. So
the function $\Phi^\xi$ is increasing along this curve, and it is
strictly increasing if and only if $m$ is not a fixed point of
$(\sq\,\xi)_M$.  But it takes on the same values at $t=0$ and $t=1$,
so $m$ must be a fixed point of $(\sq\,\xi)_M$, that is, $\xi\in\frak
g_m$. Hence $gm=g\exp(\sq\,\xi)m=m$, so $g\in G_m$.

Lastly we show (\ref{stabilizer}) implies (\ref{realorbit}). If
$(G_m)\co=(G\co)_m$, the (real) dimension of the complex orbit $G\co
m$ equals twice the dimension of the compact orbit $Gm$. Since the
tangent space at $m$ to $G\co m$ is equal to
$T_m(Gm)+J\bigl(T_m(Gm)\bigr)$, the intersection $T_m(Gm)\cap
J\bigl(T_m(Gm)\bigr)$ has to be 0, that is, $Gm$ is totally real.
\end{pf}

\begin{remark}\label{remark:connected}
The converse of the implications in the proposition are wrong. See
\cite{le:re} for a simple counterexample. If $G$ is not connected,
then the proof of the proposition shows the following implications
hold: (\ref{fixed}) $\Rightarrow$ (\ref{isotropic}) $\Rightarrow$
(\ref{stabilizer}) $\Rightarrow$ (\ref{realorbit}). Moreover, if $Gm$
is isotropic, then $\Phi(m)$ is fixed under $G^0$, where $G^0$ denotes
the component of the identity of $G$.
\end{remark}

Let $V$ be the orthogonal complement of $T_m(G\co m)$. Then using the
notation of the proposition we have an $H$-invariant orthogonal direct
sum decomposition of the tangent space:
\begin{equation}\label{equation:sum}
T_mM=\frak n\oplus J\frak n\oplus T_\mu(G\mu)\oplus V.
\end{equation}
This decomposition is symplectic in the sense that the summands $\frak
n\oplus J\frak n$, $T_\mu(G\mu)$ and $V$ are symplectic subspaces,
but it is not a complex-linear decomposition, since $T_\mu(G\mu)$
need not be $J$-invariant.

\subsection{Interpolation of K\"ahler metrics near totally real
submanifolds}\label{subsection:totallyreal}

Let $M$ be a complex manifold and let $\sigma$ be a real-valued
$C^\infty$ closed $(1,1)$-form on $M$. On any sufficiently small open
subset $O$ of $M$ we can find a potential for $\sigma$, that is, a
smooth real-valued function $u$ defined on $O$ such that $\sigma=
\sq\,\partial\bar\partial u$. Now let $X$ be any real-analytic
totally real submanifold of $M$.

\begin{theorem}\label{theorem:potential}
If the form $\sigma$ is exact in a neighbourhood of $X$\rom, then
there exists a potential for $\sigma$ defined in a \rom(possibly
smaller\rom) neighbourhood of $X$.

Assume $\sigma$ vanishes to $m$-th order on $X$\rom, that is\rom, in
local coordinates its coefficients vanish to $m$-th order on $X$\rom.
Then $\sigma$ is exact near $X$\rom, and there exists a potential for
$\sigma$ defined near $X$ which vanishes to $(m+2)$-nd order on $X$.
\end{theorem}

\begin{pf}
Because $X$ is totally real, it has a basis of Stein tubular
neighbourhoods in $M$. (See Grauert \cite{gr:on} and Reese and Wells
\cite[Theorem 2.2]{ha:ho}.) Without loss of generality we may replace
$M$ by one of these tubular neighbourhoods. Let $\alpha$ be a solution
to the equation
\begin{equation}\label{equation:poincare}
d\alpha=\sigma,
\end{equation}
and let $\beta=\alpha^{01}$ be the $(0,1)$-part of $\alpha$. It is
evident from the fact that $\sigma$ is of bidegree $(1,1)$ that
$\bar\partial\beta=0$. So we can solve the equation
\begin{equation}\label{equation:dolbeault}
\sq\,\bar\partial f=\beta,
\end{equation}
since $M$ is a Stein manifold. It is easy to check that the function
$u=f+\bar f=2\Re f$ satisfies $\sigma=\sq\,\partial\bar\partial u$.
This proves the first statement.

Now suppose $\sigma$ vanishes to $m$-th order on $X$. Then, evidently,
the restriction of $\sigma$ to $X$ is zero. Since $\sigma$ is also
closed, it follows easily from De Rham's Theorem that it is exact on
the tubular neighbourhood $M$. Let $C\colon M\times[0,1]\to M$ be the
homotopy defined by $C(p,t)=tp$, which contracts the bundle $M$ to the
zero section along the fibres. There exists a very special solution to
(\ref{equation:poincare}), namely the form
$\alpha=\int_{[0,1]}C^*\sigma$. It is not hard to check that this form
vanishes to $(m+1)$-st order on $X$. Therefore its $(0,1)$-part
$\beta=\alpha^{01}$ also vanishes to $(m+1)$-st order on $X$.  We now
want to solve (\ref{equation:dolbeault}) augmenting the order of
vanishing by one. We do this in three steps. First, we solve the
problem locally and formally. That is, we assume $X$ is an open set in
$\bold R^k$ and $M$ a strictly pseudoconvex open neighbourhood of $X$
in $\bold C^n$ (where $n\geq k$), and find a function $g$ vanishing to
$(m+2)$-nd order on $X$ such that the $(0,1)$-form
$\beta'=\beta-\sq\,\bar\partial g$ also vanishes to order $m+2$ on
$X$. Secondly, and this is the crucial point, we use H\"ormander's
$L^2$-estimates for the Cauchy-Riemann operator \cite{ho:in} to show
that locally there exists a smooth solution $g'$ to the problem
$\sq\,\bar\partial g'=\beta'$ which vanishes to $(m+2)$-nd order on
$X$. Then the function $f=g+g'$ satisfies (\ref{equation:dolbeault})
and vanishes to $(m+2)$-nd order on $X$. Thirdly, we show that the
local solutions to the problem can be glued together to obtain a
global solution, which amounts to solving a Cousin type problem.

{\it Step} 1. A {\em complexification\/} $(X\co,i)$ of the
real-analytic manifold $X$ is a complex manifold $X\co$ together with
a real-analytic map $i\colon X\to X\co$ such that for every complex
manifold $V$ and every real-analytic map $j\colon X\to V$ there
exists, for a sufficiently small open neighbourhood $O$ of $i(X)$, a
unique complex-analytic map $j\co\colon O\to V$ with $j\co\circ i=j$.
The uniqueness of the complexification (more precisely, the uniqueness
of the germ of $X\co$ at $i(X)$) is immediate from the definition; the
existence was proven by Bruhat and Whitney \cite{wh:qu}. The map $i$
is actually a closed embedding. If $j\colon X\to V$ is an embedding
and the image $j(X)$ is totally real, the complexified map $j\co$ is
an embedding (near $i(X)$).

So after shrinking the tube $M$ if necessary, we may assume that we
have inclusions $X\subset X\co\subset M$. About every point of $X$ we
can find an open neighbourhood that can be biholomorphically
identified with a strictly pseudoconvex bounded open subset $U$ of
$\bold C^k\times\bold C^l$, in such a manner that $U\cap X$ is given
by the equations $w=y=0$ and $U\cap X\co$ by $w=0$.  Here we write a
point in $\bold C^k\times\bold C^l$ as a pair $(z,w)$ with
$z=x+\sq\,y\in\bold C^k$ and $w\in\bold C^l$. We shall call a
neighbourhood with such a coordinate system a {\em distinguished\/}
neighbourhood.

Write $\beta=\sum_{\lambda=1}^{k}\zeta^\lambda\,d\bar z_\lambda +
\sum_{\lambda=1}^{l}\eta^\lambda\,d\bar w_\lambda$ and consider the
Taylor expansions at $X$ of the components $\zeta^\lambda$ and
$\eta^\lambda$:
\begin{align*}
\zeta^\lambda(x,y,w,\bar w)&\sim \sum_{|I|+|J|+|K|\geq m+2}
\zeta^\lambda_{I,J,K}(x)y^Iw^J\bar w^K,\\
\eta^\lambda(x,y,w,\bar w)&\sim \sum_{|I|+|J|+|K|\geq m+2}
\eta^\lambda_{I,J,K}(x)y^Iw^J\bar w^K,
\end{align*}
with coefficients $\zeta^\lambda_{I,J,K}$ and $\eta^\lambda_{I,J,K}$
in $C^\infty(X,\bold C)$. Here $I$, $J$ and $K$ are multi-indices
and $|I|$ denotes the norm $\sum_\lambda i_\lambda$ of
$I=(i_1,\dots,i_k)$. The fact that $\bar\partial\beta=0$ amounts to:
$$
\frac{\partial\zeta^\lambda}{\partial\bar z_\mu}=
\frac{\partial\zeta^\mu}{\partial\bar z_\lambda},\qquad
\frac{\partial\eta^\lambda}{\partial\bar z_\mu}=
\frac{\partial\zeta^\mu}{\partial\bar w_\lambda},\qquad
\frac{\partial\eta^\lambda}{\partial\bar w_\mu}=
\frac{\partial\eta^\mu}{\partial\bar w_\lambda}.
$$
Plugging the Taylor expansions of $\zeta^\lambda$ and $\eta^\lambda$
into this system of equations and inspecting the lowest-order terms in
the resulting equalities yields the following identities:
\begin{align}\label{equation:compatibility}
(i_\mu+1)\zeta^\lambda_{I+e_\mu,J,K}&=
(i_\lambda+1)\zeta^\mu_{I+e_\lambda,J,K},\notag\\
\sq\,(i_\mu+1)\eta^\lambda_{I+e_\mu,J,K}&=
(k_\lambda+1)\zeta^\mu_{I,J,K+e_\lambda},\\
(k_\mu+1)\eta^\lambda_{I,J,K+e_\mu}&=
(k_\lambda+1)\eta^\mu_{I,J,K+e_\lambda},\notag
\end{align}
for all $I$, $J$ and $K$ such that the total degree $|I|+|J|+|K|$
equals $m+1$. Here $e_\lambda$ denotes the multi-index whose entries
are all 0, except the $\lambda$-th, which is 1. (For higher-order
terms there are similar identities, but we will not need them.)

Our object is to find a smooth function $g$ such that
$\sq\,\bar\partial g=\beta$ up to terms of total degree $\geq m+3$ in
$y$, $w$ and $\bar w$. Upon substitution of the Taylor expansion of
$g$,
$$
g(x,y,w,\bar w)\sim\sum_{|I|+|J|+|K|\geq m+3}
g_{I,J,K}(x)y^Iw^J\bar w^K,
$$
we see this amounts to solving the equations
\begin{equation}\label{equation:relations}
g_{I+e_\lambda,J,K}=-\frac{2\zeta^\lambda_{I,J,K}}{i_\lambda+1},\qquad
g_{I,J,K+e_\lambda}=\frac{2\eta^\lambda_{I,J,K}}{\sq(k_\lambda+1)},
\end{equation}
for all $I$, $J$ and $K$ such that $|I|+|J|+|K|=m+2$. There are no
conditions on the terms of degree $>m+3$ in the expansion of $g$. It
is a straightforward exercise to check that the equations
(\ref{equation:relations}) with coefficients subject to the
compatibility relations (\ref{equation:compatibility}) admit solutions
$g_{I,J,K}$, where $|I|+|J|+|K|=m+3$. So if we put $g(x,y,w,\bar
w)=\sum_{|I|+|J|+|K|=m+3} g_{I,J,K}(x)y^Iw^J\bar w^K$, then $g$ is a
smooth function defined on $U$ and vanishing to $(m+2)$-nd order on
$X$, and the $(0,1)$-form $\beta'=\beta-\sq\,\bar\partial g$ also
vanishes to order $m+2$ on $X$.

{\it Step} 2. Obviously the form $\beta'$ is $\bar\partial$-closed. We
now want to find a smooth solution $g'$ defined on $U$ to the problem
$\sq\,\bar\partial g'=\beta'$ together with an order estimate. But
first note that {\em every\/} locally square integrable solution to
this equation is actually smooth. This follows from the fact that
there exists a smooth solution (\cite[Corollary 4.2.6]{ho:in}), and
that the difference of any two solutions is a $\bar\partial$-closed
function, therefore harmonic, and therefore smooth by the ellipticity
of the Laplacian on $\bold C^n$.

Let $\rho$ be the distance squared to the submanifold $X$,
$\rho(z,w)=\|y\|^2+\|w\|^2$. Because $\beta'$ vanishes to order $m+2$
on $X$, the integral $\int_U|\beta'|^2\rho^{-r}\,dzd\bar zdwd\bar w$
is finite for all $r<k+2l+m+3$. In H\"ormander's parlance $\beta'$ is
an element of the weighted $L^2$-space $L^2_{(0,1)}(U,\phi)$ with
weight $\phi=r\log\rho$. It is easy to check that for every positive
$r$ the weight function $\phi$ is plurisubharmonic. Let us take
$r=k+2l+m+2$. By Theorem 4.4.2 of H\"ormander \cite{ho:in}, we can
find a solution to the equation $\sq\,\bar\partial g'=\beta'$ such
that
$$
\int_U|g'|^2e^{-\phi}\bigl(1+\|z\|^2+\|w\|^2\bigr)^{-2}\,dzd\bar
zdwd\bar w\leq\int_U|\beta'|^2e^{-\phi}\,dzd\bar zdwd\bar w<\infty.
$$
But $g'$ is smooth as noted before, so this is only possible if $g'$
vanishes to order $\geq r-k-2l=m+2$ on $X$. The function $f=g+g'$
defined on $U$ satisfies (\ref{equation:dolbeault}) and vanishes to
order $m+2$ on $X$.

{\it Step} 3. Let $\{U_i\}$ be a Stein cover of $M$. By the previous
result we can find smooth functions $f_i$ defined on $U_i$, which
vanish to $(m+2)$-nd order on $U_i\cap X$ and satisfy
$\sq\,\bar\partial f_i=\beta|_{U_i}$. Put $f_{ij}=f_j-f_i$; then
$\bar\partial f_{ij}=0$, so $f_{ij}$ is a holomorphic function on
$U_{ij}=U_i\cap U_j$. It also vanishes to $(m+2)$-nd order on
$U_{ij}\cap X$, so, by Lemma \ref{lemma:vanish} below, it has to
vanish to order $m+2$ on $U_{ij}\cap V$, where $V$ denotes the
complexification $V=X\co$ of $X$. In other words, the collection of
$f_{ij}$'s defines a \v Cech 1-cocycle with coefficients in the
coherent sheaf $\cal I_V^{m+3}$, where $\cal I_V$ denotes the ideal
sheaf of the complex submanifold $V$. Since $M$ is Stein, Cartan's
Theorem B implies this cocycle is a coboundary (cf.\ \cite[Theorem
7.4.3]{ho:in}), so there exist holomorphic functions
$g_i\in\Gamma(U_i,\cal I_V^{m+3})$ such that $f_{ij}=g_j-g_i$.
Consider the smooth functions $f_i+g_i$ defined on $U_i$. Clearly
$f_i+g_i=f_j+g_j$ on $U_{ij}$, so $f_i+g_i=f|_{U_i}$ for a global
smooth function $f$. By construction $f$ vanishes to $(m+2)$-nd order
on $X$, and because the $g_i$'s are holomorphic, we have
$\sq\,\bar\partial f=\beta$.
\end{pf}

The proof of the theorem used the following little lemma.

\begin{lemma}\label{lemma:vanish}
Suppose $f$ is a holomorphic function on $M$ vanishing to $m$-th order
on the totally real submanifold $X$. Then $f$ vanishes to $m$-th order
on the complexification $X\co\subset M$ of $X$.
\end{lemma}

\begin{pf}
We compute in a distinguished system of coordinates $(z,w)$, writing
$z=x+\sq\,y$, as in the proof of Theorem \ref{theorem:potential}.
First we prove the statement for $m=0$. So suppose $f$ vanishes on
$X$; we have to show it vanishes on $X\co$.  The assumption implies
that the partial derivatives of $f$ along $X$,
$\partial^{|I|}f/\partial x_I$, vanish identically on $X$ for all
multi-indices $I$. Since $\partial^{|I|}/\partial z_I=
\partial^{|I|}/\partial x_I$ on holomorphic functions, we conclude
that the power series of the restriction $f|_{X\co}$ of $f$ to $X\co$
is trivial at any point of $X$. By the identity principle
$f|_{X\co}=0$.

Now suppose $f$ vanishes to order $m\geq0$ on $X$. This means the
holomorphic functions $\partial^{|I|+|J|}f/\partial z_I\partial w_J$
vanish identically on $X$ for all $I$ and $J$ with $|I|+|J|\leq m$.
Then by the previous result $\partial^{|I|+|J|}f/\partial z_I\partial
w_J=0$ on $X\co$ if $|I|+|J|\leq m$, so $f$ vanishes to order $m$ on
$X\co$.
\end{pf}

The next result says one can ``interpolate'' between two K\"ahler
metrics that agree along $X$.

\begin{theorem}\label{theorem:interpolate}
Let $dS^2$ and $ds^2$ be two smooth K\"ahler metrics on $M$. Assume
that the real-analytic totally real submanifold $X$ is compact and
that $dS^2_x=ds^2_x$ for all $x$ in $X$. Then there is an open
neighbourhood $U$ of $X$ in $M$ with the following property\rom: For
all open $U_1$ with $X\subset U_1\subset U$ there exist an open subset
$U_2$ with $X\subset U_2\subset U_1$ and a smooth K\"ahler metric
$d\tilde s^2$ on $M$ such that $d\tilde s^2=dS^2$ on $U_2$ and
$d\tilde s^2=ds^2$ on $M\backslash\bar U_1$.

In the presence of a compact group $G$ of holomorphic transformations
on $M$ leaving the submanifold $X$ and the metrics $dS^2$ and $ds^2$
invariant\rom, the metric $d\tilde s^2$ may be taken to be invariant.
If the $G$-action is Hamiltonian with respect to the K\"ahler form
$-\Im ds^2$\rom, it is Hamiltonian with respect to the K\"ahler form
$-\Im d\tilde s^2$.
\end{theorem}

\begin{pf}
For the first part of the theorem we may again assume that $M$ is a
Stein tubular neighbourhood of $X$. Let $\Omega=-\Im dS^2$ and
$\omega=-\Im ds^2$ be the K\"ahler forms corresponding to $dS^2$ and
$ds^2$, and put $\sigma=\omega-\Omega$. Then $\sigma$ vanishes to
order 0 on $X$, so by Theorem \ref{theorem:potential} there exists a
smooth function $u$ vanishing to second order on $X$ such that
$\sigma=\sq\,\partial\bar\partial u$. Let $\rho$ be the square of some
distance function on the tube $M$. Then $\rho$ vanishes to first order
on $X$, so $u$ is of order $O(\rho^{3/2})$ as $\rho$ tends to zero.
Let $\chi\colon\bold R\to[0,1]$ be a smooth function with $\chi(t)=0$
for $t\leq1$ and $\chi(t)=1$ for $t\geq2$. For $\lambda>0$ define a
smooth function $\eps$ on $M$ by
$\eps(x)=\chi\bigl(\rho(x)/\lambda^2\bigr)$. Put $M_r=\{\,x\in
M:\rho(x)<r\,\}$ and define a smooth two-form $\tilde\omega$ on $M$ by
$$
\tilde\omega=
\begin{cases}
\Omega+\sq\,\partial\bar\partial\eps u & \text{on $M_{3\lambda^2}$},\\
\omega & \text{on $M-M_{2\lambda^2}$.}
\end{cases}
$$
Then on $M_{\lambda^2}$ the form $\tilde\omega$ is equal to $\Omega$.
On $M_{3\lambda^2}$ we have
$\tilde\omega-\omega=\sq\,\partial\bar\partial(\eps-1)u$. In a
distinguished neighbourhood $U$ of a point of $X$ with coordinates
$v=(z,w)$ we can write $\partial\bar\partial(\eps-1)u=
\sum_{\alpha,\beta=1}^nf_{\alpha\beta}\,dv_\alpha\wedge d\bar v_\beta$
with
$$
f_{\alpha\beta}(v) = \frac{\partial^2}{\partial
v_\alpha\partial\bar v_\beta}\Bigl(\bigl(\eps(v)-1\bigr)u(v)\Bigr).
$$
By carrying out the differentiation one can check in a straightforward
manner that for every compact subset $K$ of $U$ the supremum of
$|f_{\alpha\beta}(v)|$ over all $v\in K\cap M_{3\lambda^2}$ is of
order $O(\lambda)$ as $\lambda$ tends to zero.  For instance, one of
the terms involved in $f_{\alpha\beta}$ is
$$
\frac{\chi''(\rho/\lambda^2)}{\lambda^4}\frac{\partial\rho}{\partial
v_\alpha}\frac{\partial\rho}{\partial\bar v_\beta}u=
\frac{\chi''(\rho/\lambda^2)}{\lambda^4}O(\rho^{5/2}),
$$
where we used $u=O(\rho^{3/2})$ and the fact that the first
derivatives of $\rho$ are of order $O(\rho^{1/2})$ as $\rho$ tends to
zero, since they vanish on $X$. Since $\chi''(\rho/\lambda^2)=0$ for
$\rho\geq2\lambda^2$, we have
$$
\sup_{\rho\leq3\lambda^2}\Bigl(
\frac{\chi''(\rho/\lambda^2)}{\lambda^4}\frac{\partial\rho}{\partial
v_\alpha}\frac{\partial\rho}{\partial\bar v_\beta}u\Bigr)=
\frac{1}{\lambda^4}O(\lambda^5)= O(\lambda).
$$
The other terms can be dealt with similarly. From the compactness of
$X$ it now follows that $\tilde\omega$ becomes arbitrarily close to
$\omega$ uniformly on $M$ as $\lambda$ tends to zero. Hence, for
$\lambda$ small enough the symmetric bilinear form
$\tilde\omega(\cdot,J\cdot)$ is positive-definite, and therefore
$\tilde\omega$ is the imaginary part of a K\"ahler metric $d\tilde
s^2$. By construction $d\tilde s^2$ is equal to $dS^2$ on
$M_{\lambda^2}$ and equal to $ds^2$ on $M-M_{2\lambda^2}$.

In the proof of the second part of the theorem we denote the Stein
tube around $X$ by $N$ to distinguish it from the whole of $M$.
Suppose the compact Lie group $G$ acts on $M$ by holomorphic
transformations leaving $X$, $dS^2$ and $ds^2$ invariant. After
averaging over $G$ we may assume the potential $u$ is invariant. It is
not hard to verify by inspecting the proof in \cite{ha:ho} that the
tube $N$ can be chosen to be invariant. If we choose an invariant
distance function, the shrunken tubes $N_r$ are also invariant. It is
clear from the definition that the form $\tilde\omega$ is then also
$G$-invariant.

Now assume that there exists a momentum map $\Phi$ for the action with
respect to the symplectic form $\omega$. By construction we have
$\tilde\omega_x=\omega_x$ for all $x$ in $X$, so by the equivariant
Darboux-Weinstein Theorem (see e.g.\ \cite[\S 22]{gu:sy}) for
sufficiently small $\lambda$ there exists a $G$-equivariant
diffeomorphism $\Gamma\colon N_{3\lambda^2}\to N_{3\lambda^2}$ fixing
the manifold $X$ such that $\Gamma^*\tilde\omega=\omega$. Then the map
$\tilde\Phi\colon N_{3\lambda^2}\to\frak g^*$ defined by
$\Gamma^*\tilde\Phi=\Phi$ is a momentum map with respect to the form
$\tilde\omega$. On $N_{3\lambda^2}-N_{2\lambda^2}$ we have
$\tilde\omega=\omega$, so there $\tilde\Phi$ differs by a locally
constant function $c$ from the $\omega$-momentum map $\Phi$. Let us
assume, as we may, that the manifolds $M$ and $X$ are connected. If
$X$ is of codimension greater than one, the subset
$N_{3\lambda^2}-N_{2\lambda^2}$, which is homeomorphic to $N-X$, is
connected. In this case $c$ is a constant, so after shifting
$\tilde\Phi$ by $c$ we can paste $\Phi$ and $\tilde\Phi$ together to
obtain a global $\tilde\omega$-momentum map for the $G$-action.

The totally real submanifold $X$ can only be of codimension one if
$\dim X=1$ and $\dim M=2$. The only Riemann surfaces $M$ that admit a
continuous group of automorphisms are $\bold P^1$, $\bold C$, $\bold
C^\times$, elliptic curves $\bold C/\Lambda$, the unit disc $\Delta$
and annuli $\Delta_r= \{\,z\in\bold C:r<|z|<1\,\}$, for $0\leq r<1$.
(Cf.\ Farkas and Kra \cite[Section V.4]{fa:ri}.) No subgroup of
$\Aut(\bold C/\Lambda)=\bold C/\Lambda$ acts on $\bold C/\Lambda$ in a
Hamiltonian fashion, so elliptic curves are out. In the other examples
the only compact connected group of automorphisms is the circle acting
in the standard way.  In each of these cases $X$ has to be a circle,
the complement of $X$ in the tube $N$ consists of two components, and
$M-X$ also consists of two components. We can therefore glue together
the two momentum maps $\Phi$ and $\tilde\Phi$ by adding appropriate
constants to $\Phi|_{M-N_{2\lambda^2}}$ on each of the two components
of $M-N_{2\lambda^2}$.
\end{pf}

\begin{remark}
If $X$ connected, then the $\tilde\omega$-momentum map $\tilde\Phi$ is
equal to $\Phi$ on $X$, and if $\Phi$ is equivariant, then so is
$\tilde\Phi$.
\end{remark}

\subsection{Holomorphic slices}\label{subsection:slices}

We now come to the main result of Section \ref{section:slices}.

\begin{theorem}[Holomorphic Slice Theorem]\label{theorem:slice}
Let $M$ be a K\"ahler manifold and let $G\co$ act holomorphically on
$M$. Assume the action of the compact real form $G$ is Hamiltonian.
Let $m$ be any point in $M$ such that the $G$-orbit through $m$ is
isotropic. Then there exists a slice at $m$ for the $G\co$-action.
\end{theorem}

If $S$ is a slice at $m$, then $gS$ is a slice at $gm$. So the theorem
implies the existence of a slice at any point $m$ such that the
$G\co$-orbit through $m$ contains an isotropic $G$-orbit. Moreover, if
$G'$ is another compact real form of $G\co$, then $G'$ is conjugate to
$G$ by some element $g$ of $G\co$, $G'=gGg^{-1}$. Then $G'$ leaves
invariant the symplectic form $g_*\omega$, and a $G'$-momentum map is
given by $\Phi'=(\Ad^*g)\circ\Phi\circ g^{-1}$, where $\Phi\colon
M\to\frak g^*$ is an $\Ad^*$-equivariant momentum map for the
$G$-action. So the choice of the compact real form is irrelevant.

\begin{trivlist}\item[\hskip\labelsep{\em Proof of Theorem
\ref{theorem:slice}}.]
We divide the proof into several steps. Using the analytic
continuation argument of Proposition \ref{proposition:orbit}, we shall
first reduce the question of the existence of a slice to the existence
of orbitally convex open neighbourhoods of the compact orbit $Gm$.
Next we consider the special case where the compact orbit $Gm$ has a
$G$-invariant neighbourhood that can be embedded in a holomorphic,
$G$-equivariant and isometric fashion into a unitary representation
space of $G$. The last step of the proof consists in showing that an
arbitrary metric $ds^2$ can always be deformed to a metric which close
to $Gm$ is the pullback of a flat metric via some embedding into a
Euclidean space, and which is still compatible with all the relevant
data. The details are as follows.

By Remark \ref{remark:connected} the vector $\Phi(m)$ is
$\Ad^*G^0$-fixed.  After shifting the momentum map we may assume that
$\Phi(m)=0$. (If $\Phi(m)$ is not fixed under the whole of $G$, then
the shifted momentum map $\Phi-\Phi(m)$ is merely $G^0$-equivariant.
This will however be sufficient in what follows.)  The tangent space
$T_mM$ at $m$ is a Hermitian vector space, which we shall identify
with standard $\bold C^n$. Then the value of the K\"ahler form
$\omega$ at $M$ is the standard symplectic form $\Omega$ on $\bold
C^n$. Let $H$ be the stabilizer of $m$ with respect to the $G$-action.
Then by Remark \ref{remark:connected} the stabilizer with respect to
the $G\co$-action is the complexification $H\co$ of $H$. The tangent
action of $H\co$ defines a linear representation $H\co\to\GL(n,\bold
C)$, the restriction of which to $H$ is a unitary representation
$H\to\U(n)$.  Let $\phi\colon O\to M$ be a local holomorphic
coordinate system on $M$ with $\phi(0)=m$ and $d\phi_0=\id_{\bold
C^n}$, where $O$ is a small $H$-invariant open ball about the origin
in $\bold C^n$. Then the pullback of the form $\omega$ is equal to
$\Omega$ at the origin.  Let $O'=\phi(O)$ and let $\psi\colon O'\to O$
be the inverse of $\phi$.  After averaging over $H$ and shrinking $O$
if necessary we may assume that $\psi$ and hence $\phi$ are
$H$-equivariant.

The tangent space to the complex orbit $G\co m$ is a Hermitian
subspace of $T_mM\cong\bold C^n$. Denote its orthogonal complement by
$V$; then $V$ is an $H\co$-invariant subspace, which can be identified
with $\bold C^l$ for some $l\leq n$. Now let $B$ be the intersection
of the ball $O$ with $V$, and let $B'$ be the image of $B$ under
$\phi$, $B'=\phi(B)$. We claim that if $B'$ is sufficiently small the
$H\co$-saturation $S'=H\co B'$ of $B'$ is a slice at $m$. (In
Snow's terminology \cite{sn:re} $B'$ is a {\em local\/} slice.) To
verify this claim we have to show that the natural map from the
associated bundle $G\co\times_{H\co}S'$ into $M$ is biholomorphic onto
an open subset of $M$. We shall show this indirectly by proving that
the map $\phi\colon B\to B'$ can be uniquely extended to a
$G\co$-equivariant map from $G\co\times_{H\co}S$ into $M$, which is
biholomorphic onto an open subset. Here $S$ is defined to be the open
subset $H\co B$ of $V$.

Let us define $E$ to be the associated bundle
$$
E=G\co\times_{H\co}V,
$$
and let $e$ be the ``base point'' $[1,0]\in E$. Consider the
$G$-equivariant map $G\times\frak m\to G\co/H\co$ sending a pair
$(g,\sq\,\xi)$ to $g\exp(\sq\,\xi)H\co$. This map descends to a
$G$-equivariant map
\begin{equation}\label{equation:mostow}
G\times_H\sq\,\frak m\to G\co/H\co,
\end{equation}
which by a refinement of the Cartan decomposition due to Mostow
\cite{mo:on1,mo:so,mo:on2} is a {\em diffeomorphism}. In other words,
every element of $G\co$ can be written as a product $g\exp(\sq\,\xi)h$
with $g\in G$, $\xi\in\frak m$ and $h\in H\co$; and if
$g\exp(\sq\,\xi)h=g'\exp(\sq\,\xi')h'$, then $g'=gk^{-1}$, $\xi'=(\Ad
k)\xi$ and $h'=kh$ for some $k\in H$. It follows that the map
$$
G\times_H(\sq\,\frak m\times V)\to G\co\times_{H\co}V
$$
sending the equivalence class $[g,\sq\,\xi,v]$ to the equivalence
class $[g\exp(\sq\,\xi),v]$ is likewise a diffeomorphism.
We conclude that the sets $U=G\exp(\sq\,D)B\simeq
G\times_H(\sq\,D\times B)$, for $D$, resp.\ $B$, ranging over all
balls about the origin in $\frak m$, resp.\ $V$, form a basis of
neighbourhoods of the compact orbit $G e$ inside the space $E$.
Furthermore, we can extend the $H$-equivariant holomorphic map
$\phi\colon B\to M$ to a $G$-equivariant holomorphic map $U\to M$ by
defining
\begin{equation}\label{equation:openset}
[g,\sq\,\xi,v]\mapsto g\exp(\sq\,\xi)\phi(v),
\end{equation}
for $g\in G$, $\xi\in D$ and $v\in B$. For simplicity we still call
this map $\phi$.  From the decomposition (\ref{equation:sum}), where
now $\mu=0$ and $\frak n=\frak m$, it is clear that
$G\times_H(\sq\,\frak m\times V)$ is nothing but the normal bundle to
the compact orbit $G m\cong G/H$ in $M$. Consequently, for $D$ and $B$
small enough $\phi\colon U\to M$ is a biholomorphic map onto an open
neighbourhood of $G m$ in $M$.  Clearly $G\co U=G\co\times_{H\co}H\co
B=G\co\times_{H\co}S$. We will prove:

\begin{claim}\label{claim:convex}
\begin{enumerate}
\item\label{model} The compact orbit $G e\subset E$ possesses a
basis of orbitally convex open neighbourhoods\rom; and
\item\label{space} The compact orbit $G m\subset M$ possesses a
basis of orbitally convex open neighbourhoods.
\end{enumerate}
\end{claim}
In view of Proposition \ref{proposition:orbit} this will imply there
is a $G\co$-equivariant biholomorphic map $\phi\co\colon
G\co\times_{H\co}S\to G\co S$ extending the map $\phi$, which will
conclude the proof of Theorem \ref{theorem:slice}. Heinzner
\cite{he:ge} gives a proof of (\ref{model}).  We shall present an
adapted version of his argument, which can be utilized to give a proof
of (\ref{space}). The argument bears a certain similarity to a
convexity argument of Kempf and Ness \cite{ke:le}. Let us start with
the simple case where $e$ is a fixed point. Then $G=H$, the space $E$
is just $\bold C^n$ and $e$ is the origin. On $\bold C^n$ we have the
constant K\"ahler metric denoted by $dS^2$, the standard symplectic
form $\Omega$ and the quadratic momentum map $\Phi_{\bold
C^n}\colon\bold C^n\to\frak g^*$ given by (\ref{equation:quadratic}).
As above, $B$ is an open ball about the origin in $\bold C^n$. By
$r(v)$ we denote the Riemannian distance of $v\in B$ to the origin,
and by $\langle\cdot,\cdot\rangle$ the positive-definite inner product
$\Re dS^2$.

\begin{lemma}\label{lemma:angle}
For all $\xi\in\frak g$ and $v\in B$ the momentum function
$\Phi_{\bold C^n}^\xi$ measures the inner product of the outward
pointing normal $\grad r^2$ to the metric sphere of radius $r$ about
the origin and the vector field $J\xi_{\bold C^n}=\grad\Phi_{\bold
C^n}^\xi$\rom, as follows\rom:
\begin{equation}\label{equation:angle}
\bigl\langle\grad r^2,\grad\Phi_{\bold
C^n}^\xi\bigr\rangle=4\Phi_{\bold C^n}^\xi.
\end{equation}

It follows that $B$ is orbitally convex with respect to the
$G\co$-action.
\end{lemma}

\begin{pf}
The path $\delta(t)=\exp(\sq\,t\xi)v$ is the trajectory of the
gradient of $\Phi_{\bold C^n}^\xi$ through $v$. On one hand,
$$
\frac{d}{dt}r^2\bigl(\delta(t)\bigr) = \bigl\langle\grad
r^2\bigl(\delta(t)\bigr),\delta'(t)\bigr\rangle = \bigl\langle\grad
r^2\bigl(\delta(t)\bigr),\grad\Phi_{\bold
C^n}^\xi\bigl(\delta(t)\bigr)\bigr\rangle.
$$
On the other hand,
\begin{align*}
\frac{d}{dt}r^2\bigl(\delta(t)\bigr)&=
\frac{d}{dt}\bigl\|\delta(t)\bigr\|^2
= \frac{d}{dt}\bigl\langle\delta(t),\delta(t)\bigr\rangle =
2\bigl\langle\delta'(t),\delta(t)\bigr\rangle = \\
&=
2\bigl\langle(\sq\,\xi)_{\bold
C^n}\bigl(\delta(t)\bigr),\delta(t)\bigr\rangle =
2\Omega\bigl(\xi_{\bold C^n}\bigl(\delta(t)\bigr),\delta(t)\bigr) =
4\Phi_{\bold C^n}^\xi\bigl(\delta(t)\bigr),
\end{align*}
where we have used (\ref{equation:quadratic}) and
(\ref{equation:grad}). Taking $t=0$ yields (\ref{equation:angle}).

Now (\ref{equation:angle}) implies that the curve $\delta(t)$ can only
enter $B$ at a point $p$ in the boundary $\partial B$ for which
$\Phi_{\bold C^n}^\xi(p)\leq0$ and leave it at a point $q\in\partial
B$ where $\Phi_{\bold C^n}^\xi(q)\geq0$. But $\delta(t)$ is also a
gradient curve of the function $\Phi_{\bold C^n}^\xi$ and so
$\Phi_{\bold C^n}^\xi$ is increasing along $\delta(t)$.  If
$\delta(t)$ is not constant, $\Phi_{\bold
C^n}^\xi\bigl(\delta(t)\bigr)$ is strictly increasing. Therefore, if
$\delta(t)$ leaves the ball $B$ at some point, it can never sneak back
in. Consequently $\{\,\delta(t):t\in\bold R\,\}\cap B$ is
connected. If $\delta(t)$ is constant it is trivially true that
$\{\,\delta(t):t\in\bold R\,\}\cap B$ is connected.
\end{pf}

Observe that the proof does not use that the metric is flat on all of
$\bold C^n$; it works for any K\"ahler metric that is flat in a
neighbourhood of the origin.

We shall make repeated use of the following result of Kempf and Ness
\cite{ke:le}. (Cf.\ also Procesi and Schwarz \cite{pr:in}.)

\begin{proposition}\label{proposition:kempfness}
Suppose $G$ acts unitarily on $\bold C^N$. Consider the complexified
representation $G\co\to\GL(N,\bold C)$. An orbit $\cal O$ of $G\co$ in
$\bold C^N$ is closed if and only if the restriction $r|_{\cal O}$ of
the length function $r$ has a stationary point. If $v\in\cal O$ is a
stationary point of $r|_{\cal O}$\rom, then
\begin{enumerate}
\item $r|_{\cal O}$ takes on its minimum at $v$\rom, and for all
$w\in\cal O$\rom, $r(w)=r(v)$ implies $w\in Gv$\rom;
\item $v$ is in the zero level set of the momentum map $\Phi_{\bold
C^N}$\rom;
\item $(G\co)_v=(G_v)\co$.\qed
\end{enumerate}
\end{proposition}

To jack up Lemma \ref{lemma:angle} we embed the homogeneous bundle $E$
equivariantly into a representation space.

\begin{lemma}\label{lemma:embedding}
There exists a $G\co$-equivariant\rom, holomorphic and proper
embedding $\iota$ of $E$ into a finite-dimensional representation
space $\bold C^N$ of $G\co$.

Choose any $G$-invariant Hermitian inner product on $\bold C^N$. Then
the sets $\iota^{-1}(B)$\rom, where $B$ ranges over the collection of
open balls about the origin in $\bold C^N$\rom, form a basis of
orbitally convex open neighbourhoods of the orbit $G e$ in $E$.
\end{lemma}

\begin{pf}
It is not hard to find an orthogonal representation of $G$ on $\bold
R^{N_1}$ for some $N_1$ containing a vector $w$ whose stabilizer is
exactly $G_w=H$. (See \cite{ja:di}.) Then the map assigning to $gH$
the vector $gw$ is a real-analytic $G$-equivariant embedding of the
homogeneous space $G/H$ into $\bold R^{N_1}$. Complexifying the
representation $G\to\operatorname{O}(N_1)$ yields a unitary
representation $G\to\U(N_1)$, which extends to a complex-linear
representation $G\co\to\GL(N_1,\bold C)$. Consider the inclusions
$Gw\subset\bold R^{N_1}$ and $G\co w\subset\bold C^{N_1}$.  Since the
$G$-representation on $\bold R^{N_1}$ is orthogonal, the tangent space
to the orbit $T_w(Gw)$ is a subspace of the tangent space to the
$(N_1-1)$-dimensional sphere about the origin in $\bold R^{N_1}$
containing $w$. It follows that the tangent space to the complexified
orbit $T_w(G\co w)=T_w(Gw)+JT_w(Gw)$ is a subspace of the tangent
space to the $(2N_1-1)$-dimensional sphere about the origin in $\bold
C^{N_1}$ containing $w$. In other words, $w$ is a critical point of
the function $r^2|_{G\co w}$, where $r^2$ is the distance to the
origin in $\bold C^{N_1}$. Proposition \ref{proposition:kempfness} now
implies that $(G_w)\co=(G\co)_w$, and that the orbit $G\co w$ is
closed in $\bold C^{N_1}$. We conclude that the map $\iota_1\colon
G\co/H\co\to\bold C^{N_1}$ sending $gH\co$ to $gw$ is an equivariant,
holomorphic and proper embedding.

Next we show how to find an embedding of the $H\co$-module $V$ into a
finite-dimensional $G\co$-module $\bold C^{N_2}$, that is, an
$H\co$-equivariant injective complex-linear map $\iota_2\colon
V\to\bold C^{N_2}$. Let $V_1$, $V_2,\dots$, $V_r$ be the irreducible
components of $V$. It is an easy consequence of the Peter-Weyl Theorem
that every irreducible $H$-module can be embedded $H$-equivariantly
into an irreducible $G$-module. (Consider the decomposition of the
left-regular representation $L^2(G)=\bigoplus_iW_i$ into
$G$-irreducibles. Decompose each of the $W_i$ into $H$-irreducibles,
$W_i=\bigoplus_{ij}Z_{ij}$. Let $\chi\colon H\to\bold C$ be the
character of some irreducible $H$-representation; pushing $\chi$
forward as a measure to $G$ gives a measure on $G$, and the
convolution product $f\to\chi*f$ defines a non-zero $H$-equivariant
projection operator $\pi$ in $L^2(G)$. Now $\pi|_{Z_{ij}}=\id$ or
$\pi|_{Z_{ij}}=0$ depending on whether or not $Z_{ij}$ has character
$\chi$. Since $\pi\neq0$ at least one of the $W_i$ has to contain a
$Z_{ij}$ with character $\chi$.) So we can find irreducible
$G$-modules $\bold C^{n_k}$ with $H$-equivariant injective
complex-linear maps $j_k\colon V_k\to\bold C^{n_k}$. Each of the
$j_k$'s is necessarily $H\co$-equivariant. We can take $\iota_2$ to be
the direct sum of the $j_k$'s.

It is now easy to check that the map $\iota\colon
E=G\co\times_{H\co}V\to\bold C^{N_1+N_2}$ mapping $[g,v]$ to
$gw+g\iota_2(v)$ is a $G\co$-equivariant, holomorphic and proper
embedding.

By Proposition \ref{proposition:kempfness} $Gw$ is precisely the
subset of vectors in $G\co w$ of minimal length. From the inequality
$\bigl\|\iota[g,v]\bigr\|^2= \bigl\|gw+g\iota_2(v)\bigr\|^2=
\|gw\|^2+\bigl\|g\iota_2(v)\bigr\|^2 \geq\|gw\|^2$ it is clear that
$Gw$ is also equal to the subset of vectors of minimal length in the
submanifold $\iota(E)$. Because of this and the $G$-invariance of the
metric on $\bold C^{N_1+N_2}$, any open ball $B$ about $0$ such that
$B\cap\iota(E)$ is nonempty contains the orbit $Gw$, and the sets
$B\cap\iota(E)$ are a basis of open neighbourhoods of the orbit
$Gw=G\cdot\iota(e)$. The second assertion of the lemma now follows
from Lemma \ref{lemma:angle} and Remark \ref{remark:trivial}.
\end{pf}

This proves part (\ref{model}) of Claim \ref{claim:convex}.

Now consider the $G$-equivariant holomorphic embeddings
$$
\bold C^N\overset{\iota}{\hookleftarrow}
U\overset{\phi}{\hookrightarrow}M,
$$
where $\phi$ is the map defined in (\ref{equation:openset}). Pulling
back the metric $ds^2$ on $M$ via $\phi$ we obtain a metric on $E$
defined on the neighbourhood $U$ of the compact orbit $Ge$. The proof
of Lemma \ref{lemma:embedding} allows us to deduce the following
stronger assertion.

\begin{lemma}\label{lemma:convex}
Suppose the linear embedding $\iota$ is isometric\rom, that is\rom,
$\iota^*dS^2=\phi^*ds^2$\rom, where $dS^2$ is the flat metric on
$\bold C^N$.  Then for any orbitally convex open subset of the form
$\iota^{-1}(B)$ contained in $U$ the image
$\phi\bigl(\iota^{-1}(B)\bigr)$ is orbitally convex in $M$.
\end{lemma}

\begin{pf}
Put $U'=\phi(U)$ and let $\psi\colon U'\to U$ be the inverse of
$\phi$. We have two $G$-invariant K\"ahler metrics on $U'$, namely
$ds^2$ and $\psi^*\iota^*dS^2$, with corresponding momentum maps
$\Phi$ and $\Phi'=\psi^*\iota^*\Phi_{\bold C^N}$.  By assumption
$ds^2$ is equal to $\psi^*\iota^*dS^2$. Moreover, $\Phi(m)=0$ and, by
Proposition \ref{proposition:kempfness}, $\Phi_{\bold
C^N}\bigl(\iota(e)\bigr)=0$.  This implies $\Phi(m)=\Phi'(m)$, and so
$\Phi=\Phi'$. Put $O=\phi\bigl(\iota^{-1}(B)\bigr)$ and pick any point
$x$ in $O$. Let $\gamma(t)\subset M$ be the curve $\exp(\sq\,t\xi)x$;
then $\gamma(t)$ is contained in $U'$ for small $t$.  Put
$v=\iota\psi(x)$ and $\delta(v)=\exp(\sq\,t\xi)v$. Let
$I=\{\,t\in\bold R:\gamma(t)\in O\,\}$ and let $I^0$ be the connected
component of $I$ containing $0$. Because the map $\iota\psi$ is
$G$-equivariant and holomorphic, we have $\delta(t)=
\exp(\sq\,t\xi)\iota\psi(x)=\iota\psi\bigl(\exp(\sq\,t\xi)x\bigr)=
\iota\psi\bigl(\gamma(t)\bigr)$ for all $t\in I^0$. (See Remark
\ref{remark:local}.) It now follows from the proof of Lemma
\ref{lemma:angle} that the curve $\gamma(t)$ can only enter the set
$O$ at time $t_0$ if $\Phi\bigl(\gamma(t_0)\bigr)\leq0$ and leave it
at time $t_1$ if $\Phi\bigl(\gamma(t_1)\bigr)\geq0$.  Because the
function $\Phi^\xi$ is increasing along $\gamma$, it follows that
$I=I^0$, i.e.\ $\gamma$ intersects $O$ in a connected set.
\end{pf}

Of course, in general the map $\iota$ will not be an isometry for the
given metric $ds^2$ on $M$. We claim, however, that we can arrange for
$\iota$ to be an isometry along the compact orbit $G e$.

\begin{lemma}\label{lemma:isometry}
The representation $G\co\to\GL(N,\bold C)$ and the $G$-invariant
Hermitian inner product on $\bold C^N$ in Lemma \ref{lemma:embedding}
can be chosen in such a way that the embedding $\iota$ is a K\"ahler
isometry at all points of the orbit $Ge$\rom, that is\rom,
$\iota^*dS^2=\phi^*ds^2$ on $T_xE$ for all $x\in G e$.
\end{lemma}

\begin{pf}
We use the notation of the proof of Lemma \ref{lemma:embedding}. Moore
\cite{mo:eq} has shown that the representation
$G\to\operatorname{O}(N_1)$ can be chosen in such a way as to make the
embedding $G/H\to\bold R^{N_1}$ an isometry of Riemannian manifolds.
(Cf.\ also \cite{mo:on}.) The associated embedding $\iota_1$ of the
complexified homogeneous space $G\co/H\co$ into $\bold C^{N_1}$ is
holomorphic and the complex structure $J$ is an orthogonal map (at
each point of $G\co/H\co$ and $\bold C^{N_1}$). So the differential
$d\iota_1$ is a unitary map $T_x(G\co/H\co)\to\bold C^{N_1}$ for all
$x\in G/H$.

We can also arrange for the embedding $\iota_2\colon V\to\bold
C^{N_2}$ to be a unitary map. Indeed, we obtained $\iota_2$ by
embedding each irreducible component $V_k$ of $V$ into an irreducible
unitary representation $\bold C^{n_k}$ of $G$. By Schur's Lemma the
invariant Hermitian inner products on $V_k$ and $\bold C^{n_k}$ are
unique up to constant multiples. By suitably rescaling the metric on
each $\bold C^{n_k}$ the embedding $\iota_2\colon\bigoplus_k
V_k\to\bigoplus_k\bold C^{n_k}$ becomes unitary.

The embedding $\iota\colon E\to\bold C^N$ is now a K\"ahler isometry
along the orbit $Ge$.
\end{pf}
With a choice of embedding as in this lemma Theorem
\ref{theorem:interpolate} tells us we can deform the metric $ds^2$ in
such a manner that $\iota$ becomes an isometry. Theorem
\ref{theorem:interpolate} plus Lemmas \ref{lemma:convex} and
\ref{lemma:isometry} therefore imply part (\ref{space}) of Claim
\ref{claim:convex}. This finishes the proof of Theorem
\ref{theorem:slice}.
\qed\end{trivlist}

Along the lines of \cite{lu:sl} one can deduce from the slice theorem
many results on the local structure of a $G\co$-action. Let us list a
few for the record.

\begin{theorem}\label{theorem:stein}
Every point in $M$ the $G$-orbit through which is isotropic
possesses a $G\co$-invariant Stein open neighbourhood.
\end{theorem}

\begin{pf}
Let $m\in M$ and suppose $Gm$ is isotropic. Put $H=G_m$. Let $S$ be a
slice at $m$ as constructed in the proof of Theorem
\ref{theorem:slice}. Then $S$ is biholomorphically equivalent to the
set $H\co B$ swept out by a ball $B$ in the tangent space $T_mS$. It
is not hard to show that $H\co B$ is Stein. We conclude $m$ has an
open $G\co$-invariant neighbourhood that is biholomorphically
equivalent to a bundle with affine base $G\co m$, Stein fibre $S$ and
reductive structure group $(G_m)\co$. By a theorem of Matsushima the
total space of this bundle is Stein.
\end{pf}

\begin{theorem}\label{theorem:subconjugate}
Let $m$ be any point in $M$ such that the $G$-orbit through $m$ is
isotropic. Then for every point $x$ nearby $m$ the stabilizer subgroup
$(G\co)_x$ is conjugate to a subgroup of $(G\co)_m$.\qed
\end{theorem}

\begin{theorem}\label{theorem:linear}
Let $m\in M$ be any fixed point of the $G\co$-action.  Then the action
of $G\co$ can be linearized in a neighbourhood of $m$ in the sense
that there exist a $G\co$-invariant open neighbourhood $U$ of $m$ in
$M$\rom, a $G\co$-invariant open neighbourhood $U'$ of the origin $0$
in the tangent space $T_mM$ and a biholomorphic $G\co$-equivariant map
$U\to U'$.
\end{theorem}

\begin{pf}
A fixed point is obviously isotropic. The result now follows
immediately from the Holomorphic Slice Theorem.
\end{pf}

\begin{remark}
This theorem was also stated by Koras \cite{ko:li}, but my proof is
different from Koras', which I have trouble understanding in places.
In particular, I fail to see a justification for his application of
the curve selection lemma.
\end{remark}

Recall that the $G\co$-action is called {\em proper\/} at the point
$m$ if for all sequences $(m_i)\subset M$ and $(g_i)\subset G$ the
following holds: If $(m_i)$ converges to $m$ and $(g_im_i)$ converges
to some point in $M$, then $(g_i)$ converges to some element of $G$.
If the action is proper at $m$, the stabilizer $(G\co)_m$ is compact.

\begin{theorem}\label{theorem:proper}
Suppose the $G$-orbit through a point $m\in M$ is isotropic.
Then the following conditions are equivalent\rom:
\begin{enumerate}
\item\label{proper} The action of $G\co$ is proper at $m$\rom;
\item\label{finite} The stabilizer $(G\co)_m$ is finite\rom;
\item\label{regular} $m$ is a regular point of the momentum map
$\Phi$.
\end{enumerate}
\end{theorem}

\begin{pf}
First we show (\ref{proper}) is equivalent to (\ref{finite}). If the
$G\co$-action is proper at $m$, the stabilizer $(G\co)_m$ is a compact
complex submanifold of $G\co$, which is a Stein manifold. Therefore
$(G\co)_m$ is finite. Conversely, assume $(G\co)_m$ is finite. Then it
is easy to see that the left action of $G\co$ on the homogeneous space
$G\co/(G\co)_m$ is proper. It follows the left $G\co$-action on the
homogeneous vector bundle $G\co\times_{(G\co)_m}V$ is also proper, $V$
being the tangent space at $m$ to a slice at $m$. By the Holomorphic
Slice Theorem the point $m$ has an invariant neighbourhood which is
equivariantly isomorphic to an invariant open subset of
$G\co\times_{(G\co)_m}V$, so the $G\co$-action on $M$ is proper at
$m$.

Next we show (\ref{finite}) is equivalent to (\ref{regular}). If
$(G\co)_m$ is finite, obviously the real stabilizer $G_m$ is also
finite, so the stabilizer subalgebra $\frak g_m$ is trivial. Now the
annihilator of $\frak g_m$ in $\frak g^*$ is equal to the range of
$d\Phi_m$ (see \cite[\S 26]{gu:sy}), so $d\Phi_m$ is surjective.
Conversely, if $d\Phi_m$ is surjective, $\frak g_m$ is trivial, so
$G_m$ is finite, so by Proposition \ref{proposition:totallyreal}
$(G_m)\co=(G\co)_m$ is finite.
\end{pf}

\begin{theorem}\label{theorem:torus}
Suppose $G$ is a torus. Then slices for the $G\co$-action exist at
all points of $M$.
\end{theorem}

\begin{pf}
If $G$ is a torus, the coadjoint representation of $G$ is trivial, so
by Proposition \ref{proposition:totallyreal} all $G$-orbits in $M$ are
isotropic. Now apply the Holomorphic Slice Theorem.
\end{pf}

Our results can also be used to give a short proof of a theorem of
Snow's \cite{sn:re}.

\begin{theorem}[Snow]\label{theorem:snow}
Let $X$ be a Stein space on which $G\co$ acts holomorphically. Let $x$
be any point in $X$ such that the orbit $G\co x$ is closed. Then there
exists a slice at $x$ for the $G\co$-action.
\end{theorem}

\begin{pf}
The first part of the proof is the same as in \cite{sn:re}. Snow
proves there exists a $G\co$-equivariant holomorphic map $h$ of $X$
into a $G\co$-representation space $\bold C^n$ that is an immersion at
$x$ (and hence at all points of the orbit $G\co x$) and whose
restriction to $G\co x$ is a proper embedding (\cite[Proposition
2.5]{sn:re}). It follows the orbit $G\co\cdot h(x)$ is closed in
$\bold C^n$, and therefore by Proposition \ref{proposition:kempfness}
the compact orbit $G\cdot h(x)$ is contained in the zero level set of
the quadratic momentum map $\Phi_{\bold C^n}$. So by Lemma
\ref{lemma:embedding} the orbit $G\cdot h(x)$ possesses a basis $\cal
U$ of orbitally convex neighbourhoods in $\bold C^n$. From the fact
that $h|_{Gx}$ is injective, that $h$ is an immersion at all points of
$Gx$ and that $Gx$ is compact, we conclude $h$ is a diffeomorphism
from a neighbourhood of $Gx$ onto a neighbourhood of $G\cdot h(x)$. It
follows that the sets $h^{-1}(U)$ for $U\in\cal U$ form a basis of
neighbourhoods of $Gx$.  By Remark \ref{remark:trivial} they are also
orbitally convex. The theorem now follows from Proposition
\ref{proposition:orbit} (or rather, the generalization of Proposition
\ref{proposition:orbit} to arbitrary complex spaces, which is just as
easy to prove; see \cite{he:ge}).
\end{pf}

\begin{remark}
One can talk of holomorphic actions and momentum maps in the setting
of K\"ahler spaces (``K\"ahler manifolds with singularities'') in the
sense of Grauert \cite[\S 3]{gr:ub}. It seems reasonable to expect
that the Holomorphic Slice Theorem can be extended to this more
general situation.
\end{remark}

\section{K\"ahler Quotients and Geometric
Quantization}\label{section:quotient}

In this section I apply the Holomorphic Slice Theorem to the study of
symplectic quotients of a K\"ahler manifold $M$. The upshot is that
such a quotient has a natural structure of an analytic space, and that
if $M$ is integral, the quotient is a complex-projective variety.

Of course, if $M$ is integral, it is a complex-projective manifold by
Kodaira's Embedding Theorem, but the embedding given by Kodaira's
theorem is usually not a symplectic embedding into projective space.
(For a simple example where it is not, consider any non-singular
$X\subset\bold CP^n$. Let $\Omega$ be the restriction of the
Fubini-Study form to $X$. For any smooth function $f$ on $X$, put
$\Omega_f=\Omega+\sq\,\partial\bar\partial f$. If $f$ is $C^2$-small,
$\Omega_f$ is a K\"ahler form. But for most $f$, for instance, those
$f$ that are not real-analytic, no holomorphic embedding of
$(X,\Omega_f)$ into any projective space $\bold CP^N$ is an isometry.)
Under the assumption that Kodaira's embedding {\em does\/} preserve
the symplectic form Kirwan \cite{ki:coh} and Ness \cite{ne:st} proved
that the symplectic quotient of $M$ agrees with a categorical quotient
of a semistable subset of $M$ in the sense of geometric invariant
theory. I show that this conclusion still holds even if Kodaira's map
is not a symplectic embedding. Thus the result says roughly that the
class of symplectic quotients of an integral K\"ahler manifold is not
bigger than the class of its algebraic quotients.  Alternatively, it
says that there are many non-equivalent symplectic structures on the
algebraic quotients of $M$.

The abovementioned result of Kirwan and Ness is a generalization of
earlier work of Guillemin and Sternberg \cite{gu:ge}, and Kempf and
Ness \cite{ke:le}. Guillemin and Sternberg dealt with the case where
the quotient of $M$ is non-singular. This case is technically simpler
mainly owing to the fact that here the action of $G\co$ is proper at
all points of the zero level set of the momentum map. (See Theorem
\ref{theorem:proper}.) Kempf and Ness handled the case of a linear
action on a Hermitian vector space. In fact, I shall reduce the
general case to that of a linear action by locally ``flattening out''
the K\"ahler metric.

Section \ref{subsection:reduction} is a discussion of quotients of
K\"ahler manifolds in the general setting of Section
\ref{section:slices}. Section \ref{subsection:integral} focuses on
the case of integral K\"ahler manifolds, placing the results of
Section \ref{subsection:reduction} in the context of geometric
invariant theory. In Section \ref{subsection:multiplicity} I rephrase
some of the results in the language of geometric quantization and show
how they lead to formul\ae\ for multiplicities of representations.

\subsection{Reduction of K\"ahler
manifolds}\label{subsection:reduction}

As in the previous section let us fix a connected K\"ahler manifold
$(M,ds^2)$ on which $G\co$ acts holomorphically and assume there
exists an equivariant momentum map $\Phi$ for the action of $G$. Let
$\lambda\in\frak g^*$.  The {\em symplectic quotient\/} or {\em
reduced \rom(phase\rom) space\/} of $M$ at the level $\lambda$ is by
definition the topological space $M_\lambda=\Phi^{-1}(G\lambda)/G$,
where $G\lambda$ is the coadjoint orbit through $\lambda$. By the
results of \cite{sj:st} $M_\lambda$ has the structure of a symplectic
stratified space. Roughly speaking, this means that $M_\lambda$ is a
disjoint union of symplectic manifolds that fit together in a nice
way, and that there is a unique open stratum, which is connected and
dense in $M_\lambda$. We want to endow $M_\lambda$ with an analytic
structure and show its stratification is analytic. We shall carry this
out only for $\lambda=0$; the general case follows from this by dint
of the ``shifting trick''. (See \cite{cu:on,sj:st}.)

Define a point $m$ in $M$ to be {\em \rom(analytically\rom)
semistable\/} if the closure of the $G\co$-orbit through $m$
intersects the zero level set $\Phi^{-1}(0)$, and denote the set of
semistable points by $M\sst$. The point $m$ is called {\em
\rom(analytically\rom) stable\/} if the closure of the $G\co$-orbit
through $m$ intersects the zero level set $\Phi^{-1}(0)$ at a point
where $d\Phi$ is surjective. The set of stable points is denoted by
$M\st$. The notions of analytic semistability and stability depend on
the K\"ahler metric and on the momentum map. If $M$ is integral, they
will turn out to be equivalent to semistability, resp.\ stability in
the sense of Mumford \cite{mu:ge} with respect to a suitable
projective embedding (Theorem \ref{theorem:gaga}).

Introduce a $G$-invariant inner product on the Lie algebra of $G$. Let
$\mu$ be the ``Yang-Mills functional'' $\|\Phi\|^2$ and let $F_t$ be
the gradient flow of the function $-\mu$. Since $\mu$ is
$G$-invariant, $F_t$ is $G$-equivariant. By Lemma 6.6 of Kirwan
\cite{ki:coh} the gradient of $\mu$ is given by
\begin{equation}\label{equation:yangmills}
\grad\mu(m)=2J\Phi(m)_{M,m},
\end{equation}
where we have identified $\Phi(m)\in\frak g^*$ with a vector in $\frak
g$ using the inner product, and where $\bigl(\Phi(m)\bigr)_{M,m}$ is
the vector field on $M$ induced by $\Phi(m)$, evaluated at the point
$m$. In particular, $\grad\mu$ is tangent to the $G\co$-orbits, so
these are preserved by the flow $F_t$. Let us call the momentum map
{\em admissible\/} if for every $m\in M$ the path of steepest descent
$F_t(m)$ through $m$ is contained in a compact set, as in Kirwan
\cite[\S 9]{ki:coh}. If $\Phi$ is admissible, the flow $F_t$ is
defined for all $t\geq0$.  Kirwan has proved $M\sst$ is the set of
points $m\in M$ with the property that the path $F_t(m)$ has a limit
point in $\Phi^{-1}(0)$. Using the ideas of Neeman \cite{ne:to} one
can show that for all $m\in M$ the limit $F_\infty(m)=
\lim_{t\to\infty}F_t(m)$ actually exists and, moreover, that the
restriction of the map $F_\infty$ to $M\sst$ is a continuous
retraction of $M\sst$ onto $\Phi^{-1}(0)$. (See also Schwarz
\cite{sc:to}.)

All proper momentum maps are admissible. Here is another simple
example, which will be important in what follows.

\begin{example}\label{example:kempfness}
Consider a linear action of $G\co$ on $\bold C^N$ with the standard
momentum map $\Phi_{\bold C^N}$ given by (\ref{equation:quadratic}).
Let $\mu_{\bold C^N}=\|\Phi_{\bold C^N}\|^2$. An easy computation
using (\ref{equation:angle}) shows that $\bigl\langle\grad
r^2,\grad\mu_{\bold C^N}\bigr\rangle=8\mu_{\bold C^N}$, where $r$
denotes the distance to the origin. Consequently, at all points of the
sphere bounding a ball $B$ about the origin the vector field
$-\grad\mu_{\bold C^N}$ points into $B$. Therefore $\Phi_{\bold C^N}$
is admissible. Let $(F_{\bold C^N})_t$ be the gradient flow of
$-\mu_{\bold C^N}$. Then it is clear that the limit map $(F_{\bold
C^N})_\infty$ retracts the set $G\co B$ onto $\Phi^{-1}_{\bold
C^N}(0)\cap B$.
\end{example}

Throughout this section we will assume $\Phi$ to be admissible. We now
collect a number of basic results on the orbit structure of $M\sst$,
which are either due to Kirwan \cite{ki:coh}, or which are refinements
of her results.

\begin{proposition}\label{proposition:semistable}
In the following \romquote closed\romunquote\ means \romquote closed
in the relative topology of $M\sst$\romunquote\ and \romquote
closure\romunquote\ means \romquote closure in $M\sst$\romunquote.
\begin{enumerate}
\item\label{open} The semistable set $M\sst$ is the smallest
$G\co$-invariant open subset of $M$ containing $\Phi^{-1}(0)$\rom, and
its complement is a complex-analytic subset\rom;
\item\label{closed} A $G\co$-orbit in $M\sst$ is closed if and only if
it intersects $\Phi^{-1}(0)$\rom;
\item\label{intersection} The intersection of a closed $G\co$-orbit
with $\Phi^{-1}(0)$ consists of precisely one $G$-orbit\rom;
\item\label{limit} For every semistable point $y$ the set
$F_\infty(G\co y)\subset\Phi^{-1}(0)$ consists of precisely one
$G$-orbit\rom;
\item\label{separate} For any pair of points $x,y\in\Phi^{-1}(0)$ that
do not lie on the same $G$-orbit there exist disjoint $G\co$-invariant
open subsets $U$ and $V$ of $M$ with $x\in U$ and $y\in V$\rom;
\item\label{closure} The closure of every $G\co$-orbit in $M\sst$
contains exactly one closed $G\co$-orbit.
\end{enumerate}
\end{proposition}

\begin{pf}
See \cite[\S 4]{ki:coh} for a proof of (\ref{open}).

We now prove (\ref{closed}). If $x$ is semistable and $G\co x$ is a
closed subset of $M\sst$, then $F_\infty(x)\in G\co x$ because the
flow $F_t$ preserves the $G\co$-orbits, and also
$F_\infty(x)\in\Phi^{-1}(0)$, so $G\co x\cap\Phi^{-1}(0)$ is
non-empty. Conversely, suppose $G\co x\cap\Phi^{-1}(0)\neq\emptyset$.
Let $(y_i)$ be a sequence in $G\co x$ converging to $y\in M\sst$. We
have to show $y\in G\co x$. Clearly $F_\infty(y)\in\overline{G\co y}$
and by continuity $F_\infty(y)=\lim_{i\to\infty}F_\infty(y_i)\in G\co
x$. Therefore $G\co y$ intersects every open neighbourhood of $G\co
x$. In particular, $y$ is contained in every $G\co$-invariant open
neighbourhood of $x$. Since $G\co x$ intersects $\Phi^{-1}(0)$, the
Holomorphic Slice Theorem tells us $x$ has a $G\co$-invariant tubular
neighbourhood $U$. Evidently, $G\co x$ is a closed subset of $U$.
Since $y\in U$, this implies $y=\lim_{i\to\infty}y_i\in G\co x$.

See \cite[\S 6]{ki:coh} for a proof of (\ref{intersection}).

For the proof of (\ref{limit}), let $x=F_\infty(y)$. Pick an arbitrary
point $z\in G\co y$. We need to show $F_\infty(z)\in Gx$. If $G$ is
connected, we can find a continuous path $\gamma\colon[0,1]\to M\sst$
with $\gamma(0)=x$, $\gamma(1)=z$ and $\gamma(t)\in G\co y$ for $t>0$.
(If $G$ is not connected, we can still do this, provided we replace
$x$ by a suitable translate $gx$, where $g\in G$.)  Consider the path
$F_\infty\circ\gamma$ contained in $\Phi^{-1}(0)$, and let $I$ be the
set of all $t$ in the unit interval such that the point $x(t)$ defined
by $x(t)= F_\infty\bigl(\gamma(t)\bigr)$ is contained in $Gx$. We
claim $I$ is open in $[0,1]$.

Indeed, suppose $t\in I$. Recall that by Lemmas
\ref{lemma:embedding}--\ref{lemma:isometry} we have a slice $S$ at $x$
with the following special properties: There exists a
$G\co$-equivariant embedding $\iota$ of $U=G\co S$ into a
$G\co$-representation space $\bold C^N$, and $U$ is equal to the set
$G\co O$, where $O$ is the inverse image $\iota^{-1}(B)$ of a
Euclidean ball about the origin in $\bold C^N$. In order not to
overburden the notation we shall identify $U$ with its image in $\bold
C^N$. Now let $(F_{\bold C^N})_t$ be the gradient flow on $\bold C^N$
associated to the function $-\mu_{\bold C^N}$ of Example
\ref{example:kempfness}. Then $(F_{\bold
C^N})_\infty=\lim_{t\to\infty}(F_{\bold C^N})_t$ retracts $U$ onto
$\Phi^{-1}_{\bold C^N}(0)\cap O$.

By choosing $B$ sufficiently small we can arrange that $O$, and hence
$U$, are contained in $M\sst$. Also, $G\co y\subset U$, since $x$ is
in the closure of $G\co y$. Since $\gamma(s)\in G\co y$ we also have
$x(s)\in\overline{G\co y}$. By (\ref{closed}) the orbit $G\co x(s)$ is
closed in $M\sst$, and hence in $U$. Furthermore, $x(s)\in O$ for $s$
sufficiently close to $t$. Therefore, $(F_{\bold
C^N})_\infty\bigl(x(s)\bigr)\in\Phi^{-1}_{\bold C^N}(0)\cap G\co
x(s)$. It now follows from (\ref{closed}) (applied to the momentum map
$\Phi_{\bold C^N}$) that for $s$ sufficiently close to $t$ the orbit
$G\co x(s)$ is closed in $\bold C^N$. By construction each $x(s)$ is
contained in the closure of the orbit $G\co y\subset U$. But in a
$G\co$-representation space each orbit contains a {\em unique\/}
closed orbit in its closure. (See e.g.\ Luna \cite[\S 1]{lu:sl}.) We
conclude that $x(s)\in G\co x$ for all $s$ close enough to $t$. Since
$\Phi\bigl(x(s)\bigr)=0$, it follows from (\ref{intersection}) that
$x(s)\in Gx$, in other words $s\in I$. Thus we have shown $I$ is open.

Obviously, $I$ is also closed and $0\in I$. It follows $I=[0,1]$, and
therefore $F_\infty(z)= F_\infty\bigl(\gamma(1)\bigr)\in Gx$. This
finishes the proof of (\ref{limit}).

To prove (\ref{separate}), observe that (\ref{limit}) implies that for
any $G$-invariant subset $A$ of $\Phi^{-1}(0)$ the preimage
$F_\infty^{-1}(A)\subset M\sst$ is $G\co$-invariant. Now suppose
$x,y\in\Phi^{-1}(0)$ and $y\not\in Gx$. Because $G$ is compact, there
exist disjoint $G$-invariant open subsets $A$ and $B$ of
$\Phi^{-1}(0)$ with $x\in A$ and $y\in B$. Then $F_\infty^{-1}(A)$ and
$F_\infty^{-1}(B)$ are disjoint $G\co$-invariant open sets containing
$x$, resp.\ $y$.

Finally, (\ref{closure}) follows immediately from (\ref{closed}) and
(\ref{separate}).
\end{pf}

Call two semistable points $x$ and $y$ {\em related\/} if the closures
of the orbits $G\co x$ and $G\co y$ intersect. (Again, ``closure''
means ``closure in $M\sst$''.) Assertion (\ref{closure}) of Proposition
\ref{proposition:semistable} implies this relation is an equivalence
relation. Write $M\sst\qu G\co$ for the quotient space and $\Pi$ for
the quotient map $M\sst\to M\sst\qu G\co$. By (\ref{separate}) above,
the space $M\sst\qu G\co$ is Hausdorff.

\begin{theorem}\label{theorem:homeomorphism}
The inclusion $\Phi^{-1}(0)\subset M\sst$ induces a homeomorphism
$M_0\to M\sst\qu G\co$.
\end{theorem}

\begin{pf}
By  (\ref{closed}) and (\ref{intersection}) of Proposition
\ref{proposition:semistable} the map $M_0\to M\sst\qu G\co$ sending a
$G$-orbit $Gm\subset\Phi^{-1}(0)$ to the equivalence class $\Pi(m)$ is
a continuous injection. By (\ref{closure}) it is a bijection.
Moreover, the inverse is induced by the retraction $F_\infty\colon
M\sst\to\Phi^{-1}(0)$ and is therefore continuous.
\end{pf}

\begin{remark}
Proposition \ref{proposition:semistable} shows that $M_0$ can also be
identified with the space of closed $G\co$-orbits in $M\sst$.
\end{remark}

Let us say that a subset $A$ of $M\sst$ is {\em saturated\/} with
respect to $\Pi$ if $\Pi^{-1}\Pi(A)=A$.  This means that for every $x$
in $A$ the closure of $G\co x$ is contained in $A$.

\begin{proposition}\label{proposition:saturated}
At every point of $\Phi^{-1}(0)$ there exists a slice $S$ such that
the set $G\co S$ is saturated with respect to the quotient mapping
$\Pi$.
\end{proposition}

\begin{pf}
Let $x$ be any point in $\Phi^{-1}(0)$. We use the notation of the
proof of part (\ref{limit}) of Proposition
\ref{proposition:semistable}. We shall show that, after shrinking $O$
if necessary, the set $U$ becomes $\Pi$-saturated. Choose a ball
$B'\subset B$ so small that the $G$-invariant neighbourhood
$O'=\iota^{-1}(B')$ of $Gx$ has the property that $F_\infty(O')\subset
O$. This is possible because the sets $\iota^{-1}(B)$ form a basis of
neighbourhoods of $Gx$ by Lemma \ref{lemma:embedding} and because
$F_\infty$ is the identity on $\Phi^{-1}(0)$.

Take any $y\in U'$. We claim that $\overline{G\co y}$ is a subset of
$U'$. Since $F_\infty(O')\subset O$, part (\ref{limit}) of Proposition
\ref{proposition:semistable} implies that $F_\infty(U')\subset O$,
where $U'=G\co O'$. In particular $F_\infty(y)\in O$, and so
$G\co\cdot F_\infty(y)\subset U$. Assertion (\ref{closed}) of
Proposition \ref{proposition:semistable} implies $G\co\cdot
F_\infty(y)$ is closed in $M\sst$, and hence in $U$. Now $(F_{\bold
C^N})_\infty$ maps $U$ into $O$, so $(F_{\bold
C^N})_\infty\bigl(F_\infty(y)\bigr)\in O$.  Moreover, since $G\co\cdot
F_\infty(y)$ is closed in $M\sst$, $(F_{\bold
C^N})_\infty\bigl(F_\infty(y)\bigr)$ sits in $G\co\cdot F_\infty(y)$.
Therefore $G\co\cdot F_\infty(y)$ is closed in $\bold C^N$ (by part
(\ref{closed}) of Proposition \ref{proposition:semistable} applied to
the momentum map $\Phi_{\bold C^N}$). But $(F_{\bold C^N})_\infty(y)$
is contained in $O'$, and $G\co\cdot(F_{\bold C^N})_\infty(y)$ is
closed in $\bold C^N$.  Moreover, both orbits $G\co\cdot(F_{\bold
C^N})_\infty(y)$ and $G\co\cdot F_\infty(y)$ are contained in the
closure of $G\co y$. It follows that $G\co\cdot
F_\infty(y)=G\co\cdot(F_{\bold C^N})_\infty(y)\subset U'$, and so
$F_\infty(y)\in U'$. We now conclude from part (\ref{limit}) of
Proposition
\ref{proposition:semistable} and the continuity of $F_\infty$ that
$\overline{G\co y}\subset U'$.
\end{pf}

{}From now on we'll identify the spaces $M\sst\qu G\co$ and $M_0$. We
want to furnish $M_0$ with a complex-analytic structure in such a way
that the quotient map $\Pi$ becomes holomorphic. The richest possible
such structure is obtained by declaring a function $f$ defined on an
open subset $A$ of $M_0$ to be holomorphic if the pullback of $f$ to
$\Pi^{-1}(A)\subset M\sst$ is holomorphic. Let $\cal O_{M_0}$ be the
sheaf of holomorphic functions on $M_0$. We claim this indeed defines
an analytic structure.

\begin{theorem}\label{theorem:analytic}
The ringed space $(M_0,\cal O_{M_0})$ is an analytic space.
\end{theorem}

\begin{pf}
Let $p\in M_0$ and let $m$ be a point in $\Phi^{-1}(0)$ sitting over
$p$. By the definition of $\cal O_{M_0}$ a neighbourhood of $p$ is
isomorphic as a ringed space to a quotient $U\qu G\co$, where $U$ is a
$\Pi$-saturated open set containing $m$, equipped with the sheaf of
$G\co$-invariant holomorphic functions. By Proposition
\ref{proposition:saturated} we may take $U$ to be of the form $G\co
S$, where $S$ is a slice at $m$.  Then $U$ can be identified with an
invariant open subset of the bundle $E=G\co\times_{(G\co)_m}V$, where
$V$ is the tangent space to a slice $S$ at $m$, and the quotient $U\qu
G\co$ can be identified with an open subset of $E\qu G\co=
V\qu(G\co)_m$. Now by a theorem of Luna \cite{lu:fo} every
$(G\co)_m$-invariant holomorphic function on the
$(G\co)_m$-representation space $V$ is a holomorphic function of the
invariant {\em polynomials\/} on $V$. Picking a finite number of
generators $(\sigma_1,\dots,\sigma_l)$ of the ring of invariant
polynomials we get a Hilbert map $\sigma\colon V\to\bold C^l$, sending
$v$ to $\bigl(\sigma_1(v),\dots,\sigma_l(v)\bigr)$. The Hilbert map
descends to a map $V\qu(G\co)_m\to\bold C^l$, which by Luna's theorem
is a closed embedding of the ringed space $V\qu(G\co)_m$. It follows
$U\qu G\co$ is isomorphic as a ringed space to an analytic subset of
an open subset of $\bold C^l$. Therefore $(M_0,\cal O_{M_0})$ is an
analytic space.
\end{pf}

The proof of this theorem shows that in a neighbourhood of the point
$p$ the quotient map $\Pi$ is equivalent to the quotient map $E\to
E\qu G\co$ of the non-singular affine $G\co$-variety $E$.

\begin{corollary}\label{corollary:affine}
The quotient map $\Pi$ is locally biholomorphically equivalent to an
affine map. In particular\rom, the fibres of $\Pi$ are affine
varieties.\qed
\end{corollary}

The Holomorphic Slice Theorem implies that if the stabilizer of a
point $m\in\Phi^{-1}(0)$ is finite, all $G\co$-orbits in an invariant
neighbourhood of $m$ must have the same dimension. From this
observation plus Theorem \ref{theorem:proper} and Proposition
\ref{proposition:semistable} one can easily deduce the following
result.

\begin{theorem}\label{theorem:stable}
If $x\in M$ is stable\rom, the orbit $G\co x$ is closed in $M\sst$ and
the stabilizer $(G_x)\co$ is finite. Let $Z$ be the set of
$m\in\Phi^{-1}(0)$ with the property that $d\Phi_m$ is surjective\rom;
then the stable set $M\st$ is equal to $F_\infty^{-1}(Z)$. Every fibre
of $\Pi|_{M\st}$ consists of a single orbit.

In particular\rom, if $0$ is a regular value of $\Phi$\rom, $M\sst$
coincides with $M\st$ and $M_0\cong M\st/G\co$\rom, the space of
stable orbits in $M$\rom, is a K\"ahler orbifold.\qed
\end{theorem}

Let $p$ be in $M_0$, let $x$ be a point in $\Phi^{-1}(0)$ mapping to
$p$ and let $(H)$ be a conjugacy class of closed subgroups of $G$.
Then $p$ is said to be of {\em $G$-orbit type\/} $(H)$ if the
stabilizer $G_x$ is conjugate to $H$ in $G$. In \cite{sj:st} we showed
that the set of all points of orbit type $(H)$ is a manifold carrying
a natural symplectic structure and that the decomposition of $M_0$
into orbit type manifolds is a stratification.

Now let $(L)$ be a conjugacy class of reductive subgroups of $G\co$.
We may assume $L=H\co$ for some closed subgroup $H$ of $G$. By
Proposition \ref{proposition:semistable} the fibre $\Pi^{-1}(p)$
contains a unique closed $G\co$-orbit, namely $G\co x$. Let us say $p$
is of {\em $G\co$-orbit type\/} $(H\co)$ if the stabilizer $(G\co)_x$
is conjugate to $H\co$ in $G\co$.

\begin{theorem}\label{theorem:stratification}
The stratification of $M_0$ by $G$-orbit types is identical to the
stratification by $G\co$-orbit types. Each stratum $\cal S$ is a
complex manifold and its closure is a complex-analytic subvariety of
$M_0$.  The reduced symplectic form on $\cal S$ is a K\"ahler form.
\end{theorem}

\begin{pf}
The first assertion boils down to showing that if $H$ and $K$ are two
closed subgroups of $G$, and $H\co$ and $K\co$ are conjugate in
$G\co$, then $H$ and $K$ are conjugate in $G$. To say that $H\co$ and
$K\co$ are conjugate in $G\co$ amounts to saying that there is a
$G\co$-equivariant diffeomorphism of homogeneous spaces $f\colon
G\co/H\co\to G\co/K\co$. By Mostow's decomposition
(\ref{equation:mostow}), for every closed subgroup $R$ of $G$ the
complexified homogeneous space $G\co/R\co$ is a homogeneous vector
bundle over $G/R$, so there exist a $G$-equivariant embedding
$\iota_R\colon G/R\to G\co/R\co$ and a $G$-equivariant retraction
$\rho_R\colon G\co/R\co\to G/R$. So the composite $\rho_K\circ
f\circ\iota_H$ is a $G$-equivariant map $G/H\to G/K$. Therefore $H$ is
conjugate (in $G$) to a subgroup of $K$. Switching the r\^oles of $H$
and $K$, we see that $K$ is conjugate to a subgroup of $H$. Therefore,
since $H$ and $K$ have finitely many components, $H$ is conjugate to
$K$ in $G$.

For quotients of affine $G\co$-varieties Luna proved in \cite[\S
III.2]{lu:sl} that each stratum is non-singular and that its closure
is a variety. In view of the fact that $\Pi$ is locally equivalent to
a quotient map of an affine variety this implies the second statement
of the theorem. (To be precise, Luna's stratification is not the same
as ours, but it is easy to see that they are the same up to connected
components.)

Let $\cal S$ be the stratum of orbit type $(H)$. It is well-known that
if $G$ acts freely on the zero level set $\Phi^{-1}(0)$ (which implies
$0$ is a regular value of $\Phi$) the reduced symplectic form is
K\"ahler. (See e.g.\ \cite{gu:ge}.)  Therefore, to prove that the
reduced symplectic form on $\cal S$ is K\"ahler, it suffices to show
that $\cal S$ can be obtained by carrying out reduction at a regular
level on some K\"ahler manifold with respect to some group action.

Let $N=N_G(H)$ be the normalizer of $H$ in $G$ and let $M_H$ be the
set of ``symmetry type'' $H$, that is, the collection of all points
whose stabilizer (with respect to the $G$-action) is exactly $H$. Then
$M_H$ is a complex submanifold of $M$, so it is K\"ahler. Moreover,
$M_H$ is $N$-invariant and $H$ acts trivially on it. Let $L=N/H$. By
Theorem 3.5 of \cite{sj:st} the momentum map maps $M_H$ into $\frak
l^*$ and the restriction $\Phi_H$ of $\Phi$ to $M_H$ is a momentum map
for the $N$-action on $M_H$.  Moreover, $0$ is a regular value of
$\Phi_H$ and the reduced space $\Phi_H^{-1}(0)$ and the stratum $\cal
S$ are symplectically diffeomorphic in a natural way.

Since $M_H$ is complex, we have a well-defined action of
$L\co=N\co/H\co$ on $M_H$. Let $(M_H)\sst=(M_H)\st$ denote the set of
points in $M_H$ stable with respect to the momentum map $\Phi_H$. By
Theorem \ref{theorem:homeomorphism} and Theorem
\ref{theorem:stable} we have a map
$(M_H)\st/L\co\cong\Phi_H^{-1}(0)/L\to M_0$, which is a homeomorphism
onto the image $\cal S$. To finish the proof, it suffices to show that
this map is biholomorphic onto $\cal S$. Since $\Phi(M_H)$ is a subset
of $\frak l^*$, (\ref{equation:yangmills}) implies that the flow of
$\|\Phi\|^2$ leaves $M_H$ invariant. Because also $\Phi_H^{-1}(0)=
\Phi^{-1}(0)\cap M_H$, we see that $(M_H)\st=M\sst\cap M_H$. Moreover,
we have a commutative diagram:
\begin{equation}\label{equation:stratum}
\begin{CD}
(M_H)\st @>\Pi>>(M_H)\st/L\co \\
@ViVV @VV{\bar\imath}V\\
M\sst @>\Pi_H>> M\sst\qu G\co,
\end{CD}
\end{equation}
where the inclusion $i$ is biholomorphic onto its image. From the
definition of the complex structures on $M\sst\qu G\co$ and
$(M_H)\st/L\co$ it now follows that $\bar\imath$ is biholomorphic onto
its image.
\end{pf}

\begin{remark}\label{remark:strata}
The orbit type stratification is the minimal real-analytic Whitney
stratification of $M_0$. However, it is not the minimal
complex-analytic stratification. This is obvious from the following
simple example. The {\em $(1,-1)$-resonance\/} is the $S^1$-action on
$\bold C^2$ defined by $e^{\sq\,\theta}\cdot(z_1,z_2)=
(e^{\sq\,\theta}z_1,e^{-\sq\,\theta}z_2)$. As a real-analytic space
the reduced space is isomorphic to the cone in $\bold R^3$ given by
$x_1^2 = x_2^2 + x_3^2$ and $x_1\geq 0$.  (See \cite{cu:on}.) There
are two strata: the vertex and the complement of the vertex. But from
the complex-analytic point of view the singularity at the vertex is
spurious: The ring of $\bold C^\times$-invariant polynomials is just
$\bold C[z_1z_2]$, so the quotient is simply $\bold C$.
\end{remark}

\subsection{The integral case}\label{subsection:integral}

The most important special case of the situation of the previous
section is that of a positive holomorphic line bundle over a complex
manifold $M$, that is, a holomorphic line bundle $\rho\colon L\to M$
with Hermitian fibre metric $\langle\cdot,\cdot\rangle$ and curvature
form $\Theta$ such that the real $(1,1)$-form
$\omega=-(2\pi\sq)^{-1}\Theta$ is K\"ahler. Recall that $\Theta$ is
the unique two-form on $M$ satisfying
$\rho^*\Theta=\bar\partial\partial\log r^2$, where $r\colon L\to\bold
R$ is the length function, $r(l)=\langle l,l\rangle^{1/2}$. The
K\"ahler class $[\omega]$ is the image of the Chern class $c_1(L)$ of
$L$ under the natural map $H^2(M,\bold Z)\to H^2(M,\bold R)$.

Now suppose that the compact group $G$ acts on $L$ by linear bundle
transformations that leave the Hermitian metric invariant. Then the
connection on $L$ is invariant and $\omega$ is invariant under the
induced action on the base $M$. This implies that for each $\xi$ in
$\frak g$ there exists a unique real-valued function $\Phi^\xi$ on $M$
such that the vector field $\xi_L$ is given by the following formula:
\begin{equation}\label{equation:lift}
\xi_L=\xi_{M,\text{hor}}+2\pi\Phi^\xi\nu_L.
\end{equation}
Here $\xi_{M,\text{hor}}$ is the horizontal lift of $\xi_M$ to $TL$
with respect to the Hermitian connection, and $\nu_L$ denotes the
vector field on $L$ generating the circle action defined by fibrewise
multiplication by complex numbers of length one.  It is not hard to
check that $\Phi^\xi$ is a Hamiltonian for the vector field $\xi_M$,
and therefore the action on $M$ is Hamiltonian.  The momentum map
$\Phi\colon M\to\frak g^*$ is automatically equivariant. The
infinitesimal action of $G$ on smooth sections $s$ of $L$ is given by:
$$
\xi\cdot s= -\nabla_{\xi_M}s+2\pi\sq\,\Phi^\xi\cdot s.
$$
As before, let us assume that the $G$-action on $M$ extends to a
$G\co$-action and that the momentum map $\Phi$ is admissible. Then the
$G$-actions on $L$ and on its smooth sections can both be uniquely
extended to actions of $G\co$, and the projection $\rho$ is
$G\co$-equivariant. (See \cite{gu:ge}.)

Note that a holomorphic section of $L$ defined over a $G\co$-invariant
open set is $G$-invariant if and only if it is $G\co$-invariant. Let
$\cal L$ be the sheaf of holomorphic sections of $L$ and define a
sheaf $\cal L_0$ on $M_0$, the {\em sheaf of invariant sections}, by
putting $\cal L_0(O)=\cal L(\Pi^{-1}(O))^G$. According to Roberts
\cite{ro:no}, $\cal L_0$ is a coherent $\cal O_{M_0}$-module. The
following result says that $\cal L_0$ is ``almost'' a holomorphic line
bundle over $M_0$.

\begin{proposition}\label{proposition:vbundle}
The sheaf $\cal L_0$ is \rom(the sheaf of sections of\rom) a
holomorphic $V$-line bundle over $M_0=M\sst\qu G\co$.
\end{proposition}

\begin{pf}
Let $p\in M_0$. We have to show there exist a neighbourhood $O$ of $p$
that can be written as a quotient of an analytic space $\tilde O$ by
the action of a finite group $\Gamma$ and a locally free sheaf
$\tilde{\cal L}_0$ of rank one over $\tilde O$ such that $\cal L_0|_O$
is isomorphic to the sheaf of $\Gamma$-invariant sections of
$\tilde{\cal L}_0$. Let $m$ be a point in $\Phi^{-1}(0)$ mapping to
$p$ and let $S$ be a slice at $m$ such that $U=G\co S$ is
$\Pi$-saturated. Put $O=\Pi(U)$; then $O\cong S\qu H\co$ as analytic
spaces. Here $H$ denotes the stabilizer $G_m$ of $m$. The group $H$
acts linearly on the fibre $L_m$. If $\eta\in\frak h$ and $l\in L_m$,
then (\ref{equation:lift}) implies $\eta_L\cdot l=0$, since
$\Phi(m)=0$. In other words, the identity component $H^0$ acts
trivially on $L_m$.

Now let $\cal N$ be the restriction of $\cal L$ to $S$ and let $s$ be
a holomorphic section of $\cal N$ that does not vanish at $m$. Then
the section $\int_{H^0}h\cdot s\,dh$ is holomorphic, $H\co$-invariant
and does not vanish at $m$, because $H^0$ acts trivially on $L_m$.
Hence we may assume $s$ to be $H^0$-invariant. Define $\tilde
O=S\qu(H^0)\co$ and let $\tilde{\cal L}_0$ be the sheaf of
$(H^0)\co$-invariant sections of $\cal N$. After shrinking $\tilde O$
if necessary, we may assume $\tilde s$ vanishes nowhere on $\tilde O$,
so $\tilde{\cal L}_0$ is a free sheaf of rank one on $\tilde O$.

By construction $O$ is the quotient of $\tilde O$ by the finite group
$H\co/(H^0)\co=H/H^0$, and $\cal L_0|_O$ is isomorphic to the sheaf of
$H/H^0$-invariant sections of $\tilde{\cal L}_0$. It follows $\cal L_0$
is the sheaf of sections of a holomorphic $V$-line bundle over $M_0$.
\end{pf}

In fact, the total space $L_0$ of this $V$-line bundle is simply the
quotient of $L|_{\Phi^{-1}(0)}$ by $G$. Furthermore, as an analytic
space $L_0$ can be identified with a quotient $L\sst\qu G\co$, where
$L\sst$ is by definition the restriction of $L$ to $M\sst$. The proofs
of these facts are sufficiently similar to the proofs in Section
\ref{subsection:reduction} that I can omit them.

To get a genuine line bundle on $M_0$, we have to replace $L$ by a
suitable power. The proof of Proposition \ref{proposition:vbundle}
shows that for every $m\in\Phi^{-1}(0)$ the image of $G_m\to\Aut(L_m)$
is a finite cyclic group. Let $q(m)$ be the order of this group.

\begin{proposition}\label{proposition:bundle}
Suppose $\Phi$ is proper. Let $q$ be the least common multiple of the
$q(m)$ for $m$ ranging over $\Phi^{-1}(0)$. Then
$(L^q)_0=(L^q)\sst\qu G\co$ is a line bundle over $M_0$ satisfying
$\Pi^*(L^q)_0=L^q|_{M\sst}$.
\end{proposition}

\begin{pf}
First note that since $\Phi$ is proper, its zero level set is compact
and so contains only finitely many orbit types. Furthermore,
$q(m)$ is not greater than the order of the component group
$G_m/(G_m)^0$. Therefore the integer $q$ is well-defined. It has the
property that $G_m$ acts trivially on the fibre $L^q_m$ for all
$m\in\Phi^{-1}(0)$. As in the proof of Proposition
\ref{proposition:vbundle} we conclude that at every semistable point
there exists a non-vanishing invariant local holomorphic section of
$L^q$. By means of these sections we can define local trivializations
of $(L^q)_0$, so $(L^q)_0$ is a holomorphic line bundle over $M_0$.

Moreover, the existence of these sections implies $G_x$ acts trivially
on $L^q_x$ for {\em all\/} semistable $x$. Using this one can easily
show that the commutative diagram
$$
\begin{CD}
L^q|_{M\sst} @>>> (L^q)_0 \\
@V\rho VV @VVV\\
M\sst @>\Pi >> M_0
\end{CD}
$$
is a pullback diagram. Therefore $\Pi^*(L^q)_0=L^q|_{M\sst}$.
\end{pf}

Grauert \cite{gr:ub} has defined a (holomorphic) line bundle $E$ over
a complex space to be {\em negative\/} if the zero section in $E$ has
a strictly pseudoconvex open neighbourhood. He called a bundle {\em
positive\/} if its dual is negative. To show that $(L^q)_0$ is
positive, we first need to discuss potentials for the reduced K\"ahler
structure.

If $m\in\Phi^{-1}(0)$, a potential for the K\"ahler form $\omega$ on
an open neighbourhood $U$ of $m$ is given by
$u=-(2\pi)^{-1}\log\langle s,s\rangle$, where $s$ is a local
holomorphic section of $L$ that does not vanish at $m$. If $G$ acts
freely on $\Phi^{-1}(0)$, the proof of Proposition
\ref{proposition:vbundle} shows we can find an invariant such section.
Then $u$ is a $G$-invariant smooth potential near $m$, and its
restriction to $U\cap\Phi^{-1}(0)$ pushes down to a smooth function
$u_0$ defined on $U_0=\bigl(U\cap\Phi^{-1}(0)\bigr)\big/G$.  It is
easy to see that $u_0$ is a potential for the reduced symplectic form
$\omega_0$.

If $G$ does not act freely on $\Phi^{-1}(0)$, it may not be possible
to find such a section, but we can certainly find an invariant local
holomorphic section $s$ of the $q$-th power of $L$ that does not
vanish at $m$. Then $u=-(2\pi q)^{-1}\log\langle s,s\rangle$ is a
$G$-invariant potential for $\omega$ near $m$, and as before its
restriction to $U\cap\Phi^{-1}(0)$ pushes down to a function $u_0$ on
the reduced space. Unfortunately, $u_0$ is not necessarily smooth or
even $C^1$ on $M_0$. (By a smooth function on $M_0$ we mean a function
that can be locally written as a differentiable function of the
holomorphic functions on $M_0$. This notion of smooth functions
differs from the one introduced in \cite{sj:st}.) This is clear from
the example in Remark \ref{remark:strata}, where $u(z_1,z_2)=
\bigl(|z_1|^2+|z_2|^2\bigr)\big/2$ and $u_0(w)=|w|/2$, with
$w=z_1z_2$. Nonetheless, we claim $u_0$ is strictly plurisubharmonic
in the sense of distributions. Recall that a continuous function $f$
on $M_0$ is {\em plurisubharmonic\/} if for all discs $D\subset\bold
C$ and all analytic maps $c\colon D\to M_0$ the distribution
$\Delta(c^*f)$ is non-negative, where $\Delta$ is the standard
Laplacian on $\bold C$. It is {\em strictly plurisubharmonic\/} if for
all smooth $g$ with compact support the function $f+\eps g$ is
plurisubharmonic for small $\eps$. (See Grauert and Remmert
\cite{gr:pl} and Lelong \cite[p. 46]{le:pl}.)

\begin{lemma}\label{lemma:plurisubharmonic}
The continuous function $u_0$ is strictly plurisubharmonic on the open
subset $U_0$ of the analytic space $M_0$. It is a potential for the
reduced K\"ahler structure in the following sense\rom: Let $\cal S$ be
any orbit type stratum in $M_0$ and let $\omega_{\cal S}$ be the
reduced symplectic form on $\cal S$. Then $u_0$ is smooth on
$U_0\cap\cal S$ and there it satisfies $\omega_{\cal
S}=\sq\,\partial\bar\partial u_0$.
\end{lemma}

\begin{pf}
The second statement is an immediate consequence of the observation
that the stratum $\cal S$ in $M_0$ of orbit type $(H)$ can be written
as a quotient $\cal S=\Phi_H^{-1}(0)/L=(M_H)\st/L\co$, where
$L=N_G(H)/H$ and $L\co$ acts properly and freely on $(M_H)\st$, as in
diagram (\ref{equation:stratum}). The function $u_0|_{\cal S}$ is
equal to the pushforward of $u|_{M_H}$ under the map
$\Phi_H^{-1}(0)\to\cal S$. Therefore it is smooth and satisfies
$\omega_{\cal S}=\sq\,\partial\bar\partial u_0|_{\cal S}$.

In particular, $u_0$ is strictly plurisubharmonic on every stratum of
$M_0$. To see it is strictly plurisubharmonic as a function on $M_0$,
we first consider the special case where $M=\bold C^n$ is a
$G\co$-representation with standard momentum map $\Phi=\Phi_{\bold
C^n}$ and standard flat metric with potential $u=\|z\|^2\big/2$.  We
embed the quotient $M_0=\bold C^n\qu G\co$ into $\bold C^l$ using
homogeneous invariant complex polynomials $\sigma_1,\dots$, $\sigma_l$
as in the proof of Theorem \ref{theorem:analytic}. We shall identify
$M_0$ with its image $\sigma(M_0)\subset\bold C^l$. Let
$w=(w_1,\dots,w_l)$ be coordinates on $\bold C^l$.  We claim $u_0$ is
strictly plurisubharmonic at the ``vertex'' $0\in M_0\subset\bold
C^l$. It suffices to show that for sufficiently small $\eps$ the
function $u_0-\eps\|w\|^2$ is plurisubharmonic close to the vertex.
(For simplicity we have written $\|w\|^2$ for the restriction of
$\|w\|^2$ to $M_0$.) Observe $u_0-\eps\|w\|^2$ is continuous on $M_0$,
so by the extension theorem for plurisubharmonic functions of Grauert
and Remmert \cite{gr:pl} it suffices to show the restriction of
$u_0-\eps\|w\|^2$ to the complement of a thin subset is
plurisubharmonic close to the vertex.  By Theorem
\ref{theorem:stratification} the complement of the topdimensional
stratum $\cal T$ is a thin subset of $M_0$.

We now exploit the fact that the cone $M_0$ is quasi-homogeneous in
$\bold C^l$. Consider the action $A$ of the positive real numbers on
$\bold C^n$ defined by scalar multiplication, $A_\lambda z=\lambda z$.
Let $d_1,\dots$, $d_l$ be the degrees of the homogeneous polynomials
$\sigma_1,\dots$, $\sigma_l$, and define an action $A$ of $\bold
R_{>0}$ on $\bold C^l$ by putting $A_\lambda(w_1,\dots,w_l)=
(\lambda^{d_1}w_1,\dots,\lambda^{d_l}w_l)$.  Then the Hilbert map
$\sigma\colon\bold C^n\to\bold C^l$ is equivariant,
$A_\lambda\circ\sigma=\sigma\circ A_\lambda$, and the stratum $\cal T$
is $A$-invariant. The Hermitian bilinear forms on $\cal T$
corresponding to the real $(1,1)$-forms $\sq\,\partial\bar\partial
u_0|_{\cal T}$ and $\sq\,\partial\bar\partial \|w\|^2\big|_{\cal T}$
are positive definite. Moreover, the flat metric on $\bold C^n$ is
conical, that is, $A_\lambda^*(\partial\bar\partial
u)=\lambda^2\partial\bar\partial u$, and therefore the induced metric
on the quotient is conical, $A_\lambda^*(\partial\bar\partial
u_0)=\lambda^2\partial\bar\partial u_0$. On the other hand,
$$
A_\lambda^*\bigl(\partial\bar\partial\|w\|^2\bigr)=
A_\lambda^*\Bigl(\sum_{ij}dw_i\wedge d\bar w_j\Bigr)=
\sum_{ij}\lambda^{d_i+d_j}\,dw_i\wedge d\bar w_j.
$$
Since $d_i\geq1$ for all $i$, we see that for sufficiently small
$\eps$ the bilinear form corresponding to
$\sq\,\partial\bar\partial\bigl(u_0-\eps\|w\|^2\bigr)\big|_{\cal T}$
is positive semidefinite on $\cal T\cap B$, where $B$ is a small ball
about the origin in $\bold C^l$. Consequently $u_0-\eps\|w\|^2$ is
plurisubharmonic on $M_0\cap B$ for small $\eps$, and so $u_0$ is
strictly plurisubharmonic at the vertex $0\in M_0$.

Now let $M$ be arbitrary and consider any point $m\in\Phi^{-1}(0)$.
Let $H=G_m$ and let $V=\bigl(T_m(G\co m)\bigr)^\perp$ be the tangent
space to the holomorphic slice at $m$. Then $V$ is a Hermitian vector
space and by the Holomorphic Slice Theorem we have an analytic
isomorphism $V\qu H\co\to M\sst\qu G\co$ defined near the vertex $0$
of $V\qu H\co$ and mapping $0$ to $\Pi(m)$. We now have two K\"ahler
metrics on the top stratum defined near $\Pi(m)$, namely the metric
$ds_0^2$ with potential $u_0$ induced by the metric on $M$, and the
metric $d\tilde s_0^2$ with potential $\tilde u_0$ induced by the flat
metric on $V$.  These metrics are not the same, but they are {\em
quasi-isometric\/} near $\Pi(m)$ in the sense that there is an
estimate of the type $C\Re ds_0^2\leq\Re d\tilde s_0^2\leq C^{-1}\Re
ds_0^2$ in $O\cap\cal T$, where $O$ is a neighbourhood of $\Pi(m)$ in
$M_0$.  From this and from the fact proved above that $\tilde u_0$ is
strictly plurisubharmonic at $\Pi(m)$, it follows that $u_0$ is also
strictly plurisubharmonic at $\Pi(m)$.
\end{pf}

We conclude the analytic space $M_0$ is a K\"ahler space as defined by
Grauert \cite{gr:ub}, if we extend Grauert's definition to include
local potentials that are not $C^2$.

\begin{theorem}\label{theorem:positive}
Assume $\Phi$ is proper. Let $\bold L$ be the line bundle
$(L^q)_0$\rom, where $q$ is as in Proposition
\ref{proposition:bundle}. Then $\bold L$ is positive in the sense of
Grauert. The reduced space $M_0$ is a complex-projective variety\rom,
a projective embedding being given by the Kodaira map $M_0\to\bold
P\bigl(H^0(M_0,\bold L^k)\bigr)$ for all sufficiently large $k$.
\end{theorem}

\begin{pf}
Let $\bold L^*$ be the dual of $\bold L$. We have to show the zero
section of $\bold L^*$ possesses a strictly pseudoconvex open
neighbourhood. The fibre metric $\langle\cdot,\cdot\rangle$ on $L^q$
pushes down to a fibre metric $\langle\cdot,\cdot\rangle_0$ on $\bold
L$. Let $\langle\cdot,\cdot\rangle_0^*$ be the fibre metric on $\bold
L^*$ obtained by duality. On $L^q$ we have the distance function
$r(l)=\langle l,l\rangle^{1/2}$. Let $r_0$ and $r_0^*$ be the
corresponding functions on $\bold L$, resp.\ $\bold L^*$. Let
$\Delta\subset\bold L^*$ be the tubular domain
$\{\,l:r_0^*(l)\leq1\,\}$.

In a local trivialization $(z,\zeta)$ of $\bold L$ over an open subset
$O$ of $M_0$ we can write $r_0(z,\zeta)=h_0(z)|\zeta|^2$ for a certain
positive function $h_0$ on $O$. We can use the coordinates
$\bigl(z,\bar\zeta\bigr)$ to trivialize $\bold L^*$ over $O$; then
$r_0^*\bigl(z,\bar\zeta\bigr)=h_0(z)^{-1}|\zeta|^2$. Also,
$\Delta\cap(\rho_0^*)^{-1}(O)$ is given by $|\zeta|^2\leq
h_0(z)$, where $\rho_0^*\colon\bold L^*\to M_0$ is the bundle
projection. Up to a positive constant factor the function $u_0=-\log
h_0$ is a local potential for the reduced K\"ahler structure, so by
Lemma \ref{lemma:plurisubharmonic} it is strictly plurisubharmonic. It
follows immediately that $\Delta\cap(\rho_0^*)^{-1}(O)$ is
strictly pseudoconvex in $(\rho_0^*)^{-1}(O)$. Thus we have
shown $\Delta$ is a strictly pseudoconvex subset of $\bold L^*$.

For the second part of the theorem, apply Grauert's generalization of
Kodaira's Embedding Theorem, \cite[\S 3]{gr:ub}, Satz 2.
\end{pf}

Let us call a point $x\in M$ {\em algebraically semistable\/} if there
exists an invariant global holomorphic section $s\in\Gamma(M,L^{l})^G$
of some power $L^{l}$ of $L$ such that $s(x)\neq0$. The point $x$ is
called {\em algebraically stable\/} if in addition $G\co$ acts
properly on the open set $\{\,x\in M:s(x)\neq0\,\}$. If $M$ is
algebraic, for instance, if $M$ is compact, these notions coincide
with the ones introduced by Mumford \cite{mu:ge} (except that Mumford
uses the term ``properly stable'' where most authors nowadays use
``stable'').

\begin{theorem}\label{theorem:gaga}
If $\Phi$ is proper\rom, the quotient map $\Pi\colon M\sst\to M_0$ and
the inclusion $M\sst\subset M$ induce isomorphisms
$\Gamma(M_0,L_0)\cong\Gamma(M\sst,L)^G\cong\Gamma(M,L)^G$. It follows
that a point in $M$ is analytically \rom(semi\rom)stable if and only
if it is algebraically \rom(semi\rom)stable.
\end{theorem}

\begin{pf}
The first isomorpism follows from Proposition
\ref{proposition:vbundle}. The second isomorphism follows from the
observation, essentially due to Guillemin and Sternberg \cite{gu:ge},
that the norm of an invariant holomorphic section $s$ of $L$ is
increasing along the trajectories of $-\grad\mu$. Indeed, for any
invariant holomorphic section $s$ defined on a $G\co$-invariant open
subset and any $\xi\in\frak g$ we have $J\xi_M\langle
s,s\rangle=-4\pi\Phi^\xi\langle s,s\rangle$ (see \cite{gu:ge}), so
using (\ref{equation:yangmills}) we get for any $x\in M$:

\begin{align}\label{equation:maximum}
\begin{split}
\frac{d}{dt}\bigl\langle s(F_tx),s(F_tx)\bigr\rangle&=
-\grad\mu\bigl(\langle s,s\rangle\bigr)(F_tx)=\\
&=-2J\Phi(F_tx)_M\bigl(\langle s,s\rangle\bigr)(F_tx)=\\
&=8\pi\bigl\|\Phi(F_tx)\bigr\|^2\bigl\langle
s(F_tx),s(F_tx)\bigr\rangle=\\ &=8\pi\mu(F_tx)\bigl\langle
s(F_tx),s(F_tx)\bigr\rangle\geq\\
&\geq0.
\end{split}
\end{align}
It follows that for all $x\in M\sst$ the restriction of the function
$\langle s,s\rangle$ to $\overline{G\co x}$ takes on its maximum at
$F_\infty x$.  Therefore, if $s$ is defined on all of $M\sst$,
$\langle s,s\rangle$ is bounded on $M$, since $\Phi$ is proper. An
application of Riemann's Extension Theorem now gives
$\Gamma(M\sst,L)^G\cong\Gamma(M,L)^G$.

Now suppose $x\in M$ is algebraically semistable. Then there exists
$s\in\Gamma(M,L^{l})^G$ such that $s(x)\neq0$. As $t$ tends to
infinity, $F_tx$ approaches the critical point $F_\infty x$, so
$\grad\mu(F_tx)\to0$.  Letting $t\to\infty$ in
(\ref{equation:maximum}) we obtain
$$
\mu(F_\infty x)\bigl\langle s(F_\infty x),s(F_\infty
x)\bigr\rangle=0.
$$
But since $F_\infty x\in\overline{G\co x}$, we have $s(F_\infty
x)=s(x)\neq0$.  Therefore $\mu(F_\infty x)=0$, in other words $x$ is
analytically semistable.

Conversely, suppose $x\in M$ is analytically semistable, that is,
$\Phi(m)=0$, where $m=F_\infty x$. By Proposition
\ref{proposition:bundle}, $\bold L=(L^q)_0$ is an ample bundle on
$M_0$, so $\bold L^r$ is generated by global sections for big $r$. Let
$s_0$ be a global section of $\bold L^r$ with $s_0\bigl(\Pi(x)\bigr)=
s_0\bigl(\Pi(m)\bigr)\neq0$.  By the first part of the theorem
(applied to the bundle $L^{rq}$) we can lift $s_0$ to a global
invariant section $s\in\Gamma(M,L^{rq})^G$.  Evidently, $s(x)\neq0$,
so $x$ is algebraically semistable.

One proves the equivalence of analytic and algebraic stability in a
similar way, using Theorem \ref{theorem:stable}.
\end{pf}

To summarize, the set $M-M\sst$ can be characterized as the collection
of points where all invariant global holomorphic sections of all
powers of $L$ vanish. Also, the algebraic structure of $M_0$ depends
only on the line bundle $L$ and the lift of the $G$-action to $L$, not
on the symplectic form $\omega$ or the momentum map $\Phi$. (The
symplectic structure of course does depend on $\omega$ and $\Phi$.) If
$M$ is compact, we conclude that as a projective variety $M_0$ is
nothing but the quotient defined by Mumford \cite{mu:ge}.

\subsection{Multiplicity formul\ae}\label{subsection:multiplicity}

In this section I have a stab at the ``geometric multiplicity theory''
of singular symplectic quotients. Let $M$ be a K\"ahler manifold
furnished with $G$-equivariant ``prequantum data''
$\bigl(L,\langle\cdot,\cdot\rangle\bigr)$ and momentum map $\Phi$ as
in the previous section. As before, let us assume that the $G$-action
on $M$ extends to a $G\co$-action. Let us also suppose for simplicity
that the map $\Phi$ is proper and that the group $G$ is connected.
Ideally, one would like to show that as a $G$-representation the space
of sections $\Gamma(L)$, sometimes called a ``quantization'' of $M$,
is a symplectic invariant of $M$, in other words, that it is
independent of the choice of the complex structure and the line bundle
on $M$. One way of doing this would be to express the multiplicities
of the unitary irreducible representations occurring in $\Gamma(L)$ in
terms of symplectic data involving the reduced phase spaces
$M_\lambda=\Phi^{-1}(G\lambda)/G$, where $\lambda$ ranges over the
positive weights in the dual of a maximal torus of $\frak g$.
Guillemin and Sternberg \cite{gu:ge} carried this out for those
weights $\lambda$ for which $G$ acts freely on $\Phi^{-1}(G\lambda)$.
(Then the reduced space $M_\lambda$ is non-singular.) Heckman
\cite{he:pr} had earlier obtained related results in the important
special case where $M$ is a coadjoint orbit of a big group containing
$G$ as a subgroup. Using the results of the previous section, we can
generalize their results.

If we regard the space of sections $\Gamma(M_0,L_0)$ of the $V$-line
bundle $L_0$ as the quantization of the reduced space $M_0$, then
Theorem \ref{theorem:gaga} bears out the principle that quantization
should commute with reduction.

\begin{remark}
The theory of geometric quantization is usually phrased in terms of
polarizations, that is, involutive Lagrangian subbundles of the
complexified tangent bundle of $M$. The quantization of $M$ with
respect to a polarization $P$ is then the space of polarized sections
of $L$, that is, sections $s$ such that $\nabla_Xs=0$ for all vectors
$X$ tangent to the conjugate subbundle $\bar P$. It appears to be
difficult to make sense of the notion of a polarization on a singular
space, such as a symplectic quotient. In the case of a K\"ahler
quotient $M_0$, however, the sheaf of holomorphic functions $\cal
O_{M_0}$ seems to be a workable substitute for a polarization.
\end{remark}

\begin{remark}
Properly speaking, the quantization of $M$ is not just the space of
holomorphic sections of $L$, but the virtual representation
$\bigoplus_i(-1)^iH^i(M,L)$, including all cohomology groups with
coefficients in $L$. One might wonder whether it is true that
reduction commutes with quantization in this broader sense, that is,
$\bigoplus_i(-1)^iH^i(M,L)^G=\bigoplus_i(-1)^iH^i(M_0,L_0)$. It is not
hard to prove that $H^i(M_0,L_0)$ is isomorphic to $H^i(M\sst,L)^G$
for all $i$, but I don't know if $H^i(M\sst,L)^G$ is isomorphic to
$H^i(M,L)^G$ for $i>0$.
\end{remark}

Theorem \ref{theorem:gaga} obviously implies that the dimension of
$\Gamma(M_0,L_0)$ is equal to the multiplicity of the one-dimensional
trivial representation in $\Gamma(M,L)$. Let us briefly recall from
\cite{gu:ge} how this statement generalizes to arbitrary
multiplicities by dint of the shifting trick. Choose a maximal torus
$\frak t$ in $\frak g$ and a positive Weyl chamber $\frak t^*_+$ in
$\frak g^*$.  For every positive weight $\lambda\in\frak t^*_+$ the
coadjoint orbit $G\lambda$ is a K\"ahler manifold carrying a naturally
defined Hermitian line bundle $V_\lambda$, and the Borel-Weil Theorem
asserts that $\Gamma(G\lambda,V_\lambda)$ is the unitary irreducible
representation with highest weight $\lambda$.  Let $G\lambda^-$ be the
orbit $G\lambda$ with the opposite symplectic and complex structures
and consider the K\"ahler manifold $M\times G\lambda^-$.  Let $\pi_M$
and $\pi_{G\lambda}$ denote the projections of $M\times G\lambda^-$ on
the respective factors and let $V_\lambda^*$ be the dual of
$V_\lambda$. Then the Hermitian line bundle
$\pi_M^*L\otimes\pi_{G\lambda}^*V_\lambda^*$ prequantizes $M\times
G\lambda^-$. The reduced space at 0 of $M\times G\lambda^-$ can be
identified with $M_\lambda$, the reduced space of $M$ at the orbit
$G\lambda$, and it comes equipped with a $V$-line bundle
$L_\lambda=L|_{\Phi^{-1}(G\lambda)}\bigm/G$. By Theorem
\ref{theorem:gaga}, $\Gamma(M_\lambda,L_\lambda)$ is isomorphic to the
space of $G$-invariants in $\Gamma(M\times
G\lambda^-,\pi_M^*L\otimes\pi_{G\lambda}^*V_\lambda^*)$. The K\"unneth
theorem for coherent sheaves \cite{sa:ku} now implies the following
assertion.

\begin{theorem}\label{theorem:shift}
For every positive weight $\lambda$ of $G$ the space of sections
$\Gamma(M_\lambda,L_\lambda)$ is naturally isomorphic to the space of
intertwining operators
$$
\Hom\bigl(\Gamma(G\lambda,V_\lambda),\Gamma(M,L)\bigr)^G.
\qed
$$
\end{theorem}

\begin{corollary}
If the orbit $G\lambda$ does not lie in the image of the
momentum map, the irreducible representation corresponding to
$\lambda$ does not occur in $\Gamma(M,L)$.\qed
\end{corollary}

Let us write $\mu(\lambda,L)$ for the multiplicity of the
representation with highest weight $\lambda$ occurring in
$\Gamma(M,L)$. By Theorem \ref{theorem:shift} $\mu(\lambda,L)$ is
equal to the dimension of $\Gamma(M_\lambda,L_\lambda)$. By Theorem
\ref{theorem:positive}, for certain $q$ (possibly depending on
$\lambda$) the sheaf $\bold L=(L^q)_\lambda$ is an ample line bundle
on the projective variety $X=M_\lambda$, so for all sufficiently large
$r$ we have $H^i(X,\bold L^r)=0$ for $i>0$. Then $\mu(rq\lambda,\bold
L^r)$ is equal to the Euler characteristic $\chi(X,\bold L^r)$, so by
the Hirzebruch-Riemann-Roch Theorem of Baum, Fulton and MacPherson
\cite{ba:ri1,ba:ri2}
\begin{equation}\label{equation:multiplicity}
\mu(rq\lambda,\bold L^r)=\chi(X,\bold L^r)=\eps\bigl(\ch\bold
L^r\cap\tau(X)\bigr).
\end{equation}
Here $\tau(X)$ denotes the homological Todd class of $X$, $\ch\bold
L^r$ denotes the Chern character of $\bold L^r$ and $\eps$ is the
augmentation (the map $H_\bu(X)\to\bold C$ induced by mapping $X$ to a
point). If $X$ is non-singular, (\ref{equation:multiplicity}) comes
down to the classical Hirzebruch-Riemann-Roch Theorem. As was pointed
out by Guillemin and Sternberg, $\mu(rq\lambda,\bold L^r)$ is then a
symplectic invariant, i.e., independent of the complex structure and
the line bundle $L$ on $M$. It seems likely that this is also true if
$X$ is singular. I cannot quite prove this, but here follows some
evidence.

For arbitrary singular spaces, the Todd class appears to be
intractable, but for spaces with quotient singularities, such as $X$,
the situation is simpler. Namely, a theorem of Boutot \cite{bo:si}
asserts that quotient singularities are {\em rational}, i.e.,
\begin{equation}\label{equation:rational}
f_*\cal O_Y=\cal O_X \qquad\text{and}\qquad R^if_*\cal
O_Y=0\quad\text{for }i>0,
\end{equation}
where $f\colon Y\to X$ is a resolution of singularities of $X$.
Therefore, by the functorial properties of the Todd class,
$\tau(X)=f_*\bigl(\tau(Y)\bigr)$, so by
(\ref{equation:multiplicity}) $\mu(rq\lambda,\bold L^r)=\chi(X,\bold
L^r)$ is equal to $\chi(Y,f^*\bold L^r)$.

In general, it is difficult to write down a desingularization of $X$,
but Kirwan \cite{ki:pa} has explicitly constructed a ``partial''
resolution $p\colon\tilde X\to X$. It has all the properties of a
desingularization, except that it is not a smooth projective variety,
but a complex-projective $V$-manifold (or orbifold). It is easy to see
that the vanishing property (\ref{equation:rational}) also holds for
the partial resolution $p\colon\tilde X\to X$, and therefore
$\mu(rq\lambda,\bold L^r)=\chi(\tilde X,p^*\bold L^r)$. To construct a
partial resolution, one performs a certain sequence of blowups on
$M\times G\lambda^-$ at $G\co$-invariant submanifolds, yielding a
projective manifold $\tilde M$ with a $G\co$-action. The symplectic
form on $M\times G\lambda^-$ pulls back to a degenerate $(1,1)$-form
$\tilde\omega$ on $\tilde M$, which descends to a degenerate
$(1,1)$-form $\tilde\omega_\lambda$ on $\tilde X$. The class of
$q\tilde\omega_\lambda$ is the Chern class of the pull-back $p^*\bold
L$ of $\bold L$. To get a K\"ahler form on $\tilde M$, one adds to the
pullback of $\omega$ at each stage in the sequence of blowups a small
$(1,1)$-form $\sigma_\eps$ supported on a neighbourhood of the
exceptional divisor, such that the class of $\sigma_\eps$ is equal to
$\eps$ times the dual class of the exceptional divisor. One then
obtains $\tilde X$ by taking the quotient of $\tilde M$; $\tilde
X=\tilde M\sst\qu G\co\simeq\tilde M_0$. See \cite{ki:pa} for the
details. Using Kawasaki's formula \cite{ka:ri}, we can now write the
multiplicity as a sum of integrals (still assuming that $r$ is
sufficiently large):
\begin{equation}\label{equation:integral}
\mu(rq\lambda,\bold L^r)=\int_{\tilde X}
e^{rq\tilde\omega_\lambda}\wedge
\det\frac{\sq\,R/2\pi}{1-e^{-\sq\,R/2\pi}} + \Sigma,
\end{equation}
where $R$ denotes the curvature two-form of $\tilde X$ with respect to
the K\"ahler metric, and where $\Sigma$ denotes a sum of contributions
from the singular strata in $\tilde X$. In a local $V$-manifold chart,
$\Sigma$ can be written as a sum of integrals over fixed-point
manifolds, as in the holomorphic Lefschetz formula. Since the
right-hand side is a function of the cohomology class of
$\tilde\omega_\lambda$ only, we conclude that it is a polynomial
function of $r$, and hence $\mu(rq\lambda,\bold L^r)$ is polynomial in
$r$ for $r$ large. If the dimension of $X$ is $2n$, the form
$\tilde\omega_\lambda$ enters in the term $\Sigma$ with exponents less
than $n$. The highest-order term in the multiplicity is therefore the
term in $r^n$, and the coefficient is $q^n\int_{\tilde
X}\tilde\omega_\lambda^n/n!$, which is equal to $q^n\vol X$, where
$\vol X$ is the symplectic volume of the top-dimensional stratum of
$X=M_\lambda$. In particular, we see that the highest-order term is a
symplectic invariant.

The following points seem to call for further clarification: (i) Is it
really necessary in (\ref{equation:multiplicity}) and
(\ref{equation:integral}) to replace the $V$-bundle $L_\lambda$ by the
line bundle $\bold L=(L^q)_\lambda$? In other words, does the
Riemann-Roch formula of Baum, Fulton and MacPherson work for
$V$-bundles on $M_\lambda$? It seems reasonable to guess that
$\chi(M_\lambda,L_\lambda)$ equals $\chi(\tilde
M_\lambda,f^*L_\lambda)$, which can then be computed by Kawasaki's
recipe.  (ii) Under what conditions does vanishing of the cohomology
of $M$ with coefficients in $L$ imply vanishing of the cohomology of
the quotients? (iii) More importantly, the right-hand side of
(\ref{equation:integral}) makes sense even if the symplectic manifold
$M$ does not carry a complex structure: Kirwan's partial resolution
can be defined for any singular symplectic quotient, and all one needs
to write down the form representing the Todd class is an {\em
almost}-complex structure compatible with the symplectic form. It
would be interesting to find out in how far (\ref{equation:integral})
can be generalized to this more general situation.

\ifx\undefined\bysame
\newcommand{\bysame}{\leavevmode\hbox to3em{\hrulefill}\,}
\fi

\end{document}